\documentclass[a4paper,11pt]{article}
\usepackage{amsfonts,amssymb,latexsym}
\usepackage{amsmath}
\usepackage{amsthm}
\usepackage{color}
\usepackage{graphicx}
\usepackage{subcaption}
\usepackage{epsfig,soul}
\usepackage{graphicx,graphics,float}

\newtheorem{thm}{Theorem}[section]
\newtheorem{defn}[thm]{Definition}
\newtheorem{cor}[thm]{Corollary}

\newtheorem{lem}[thm]{Lemma}

\newtheorem{rem}[thm]{Remark}

\newcommand{\be}{\mbox{$\mathbb E$}}

\newcommand{\cf}{\mbox{$\cal F$}}

\newcommand{\cL}{\mbox{$\cal L$}}

\newcommand{\br}{\mbox{$\mathbb R$}}
\newcommand{\bn}{\mbox{$\mathbb N$}}
\newcommand{\bp}{\mbox{$\mathbb P$}}

\newcommand{\de}{\delta}

\newcommand{\si}{\sigma}

\newcommand{\Om}{\Omega}

\newcommand{\R}{{\mathbb R}}
\newcommand{\Z}{{\mathbb Z}}
\newcommand{\N}{{\mathbb N}}

\newcommand{\Pb}{{\mathbb P}}

\newcommand{\calB}{{\cal B}}

\newcommand{\bfa}{\mathbf{a}}
\newcommand{\bfi}{\mathbf{i}}
\newcommand{\bfj}{\mathbf{j}}
\newcommand{\bfk}{\mathbf{k}}

\newcommand{\bfx}{\mathbf{x}}

\newcommand{\bftheta}{\boldsymbol{\theta}}

\makeatletter
\@addtoreset{equation}{section}

\makeatother
\def\enpf{
   {  \parfillskip=0pt\hfil \mbox{$\Box$} \par\bigskip  }
   }

\makeatletter

\begin{document}

\bibliographystyle{plain}

\begin{center}
{\Large{\bf A probabilistic model for interfaces  \\[.1in]
in a martensitic phase transition}}
\end{center}

\vspace*{.2 true in}
\begin{center}
{\large P. Cesana}\footnote{Institute of Mathematics for Industry, Kyushu University,
744 Motooka, Fukuoka 819-0395, Japan and Department of Mathematics and Statistics, La Trobe University, 3086 Bundoora, VIC, Australia} and {\large B.M.~Hambly}\footnote{Mathematical Institute, University
of Oxford, Radcliffe Observatory Quarter, Woodstock Road, Oxford OX2 6GG, UK.}
\end{center}
\vspace*{.1in}

\centerline{\today}
\vspace*{.1 true in}

\begin{abstract}
We analyse features of the patterns formed from a simple model for a martensitic phase transition. 
This is a fragmentation model that can be encoded by a general branching random walk. An important quantity is the distribution of the lengths of the interfaces in the pattern and we establish limit theorems for some of the asymptotics of the interface profile. We are also able to use a general 
branching process to show almost sure power law decay of the number of interfaces of at least a 
certain size. We discuss the numerical aspects of determining the behaviour of the density profile 
and power laws from simulations of the model.
\end{abstract}

\begin{small}

\end{small}

\section{Introduction}

\begin{figure}[h!]
\centering%
\includegraphics[width=7cm]{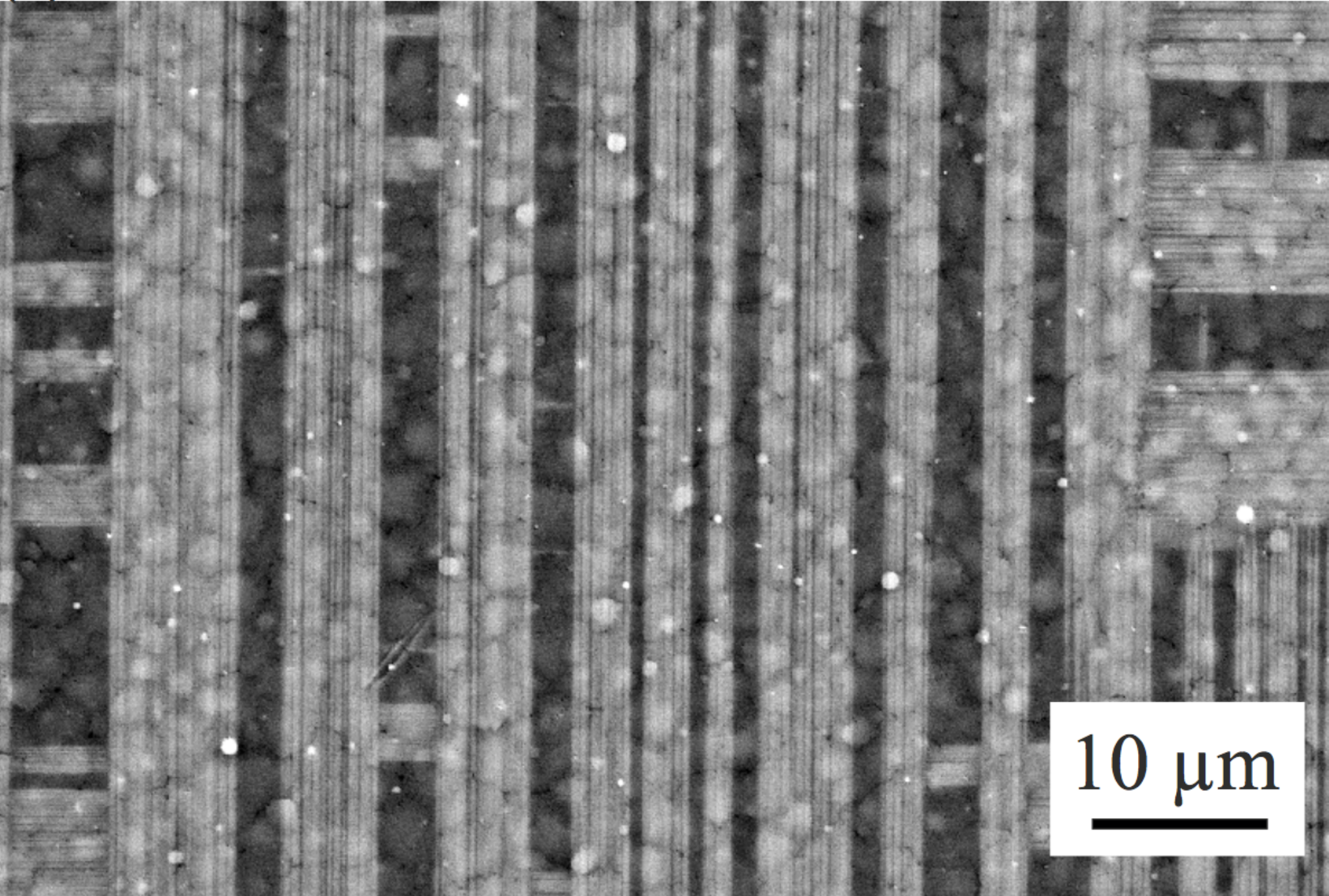}
\caption{
Snapshot from a temperature series of scanning electron micrographs of the alloy NiMnGa.
Here the dark grey area corresponds to untransformed material (austenite) and the light grey 'needles' 
evolving in vertical or horizontal direction correspond to the domain already transformed into a 
martensitic phase. Courtesy of Robert Niemann and reproduced with the permission of the author.
}\label{1609181210}
\end{figure}

In this paper, we are concerned with the stochastic modelling and analysis of the self-similar 
patterns observed in a class of elastic crystals undergoing a shape-change at the level of the
atomic lattice as a consequence of a first-order thermodynamic phase-transformation.
This is the Austenite-to-Martensite (AM) phase transformation, a 
solid-to-solid transition observed in steel and shape-memory alloys 
driven by temperature or an applied mechanical stress \cite{Bhatta}.
In a temperature-driven transformation the material evolves from austenite, the highly symmetric and 
highly homogenous crystal phase (typically the cubic lattice), as temperature
decreases via a series of ``avalanches" from one energetically metastable state to another.
%
The energy releases associated with the avalanches are manifested in the nucleation and growth of thin 
regions where the molecules are arranged locally according to a lower-symmetry lattice (martensite) as 
well as in acoustic waves generated by their moving interfaces. As a result, an intricate highly inhomogeneous
pattern soon emerges  populated by sharp interfaces that separate thin plates composed of mixtures of 
different martensitic phases (i.e., rotated copies of the low-symmetry lattice) and planes separating 
the austenite from martensite, what in the materials science literature is known as a microstructure 
(see Figure~\ref{1609181210} and also a \href{https://www.youtube.com/watch?v=Bwy8xC3hous}{movie} 
supplementary to \cite{Song}). An important point is that the Austenite-to-Martensite transformation is 
purely structural. It does not involve chemical reactions or diffusion but only a morphological 
change in the atomic arrangements. While the nucleation site of the interfaces is influenced by the state of 
disorder of the system, lattice defects and impurities which typically  requires a probabilistic modelling 
approach, their orientation and directions of propagation are determined by the kinematic compatibility of 
the various crystal phases, which are defined by rigid algebraic relations incorporated in the symmetry 
properties of the transformed lattice.





Statistical analysis of the acoustic activity captured during the transformation \cite{Salje} shows that 
several macroscopic variables such as amplitude, duration and energy of the acoustic waves exhibit a power 
law behaviour and a characteristic exponent. More experimental evidence \cite{Planes13} seems to suggest 
that the exponents measured are indeed characteristic of an entire class of materials undergoing a specific 
AM transformation, 
leading to the conclusion, which constitutes the ansatz of this work, that materials with different 
chemical composition may be grouped into universality classes identified uniquely by the symmetry 
reduction of the AM transformation.

Although of great relevance for technological applications and despite the substantial literature on the 
analysis and prediction of the statics of this transformation, it is remarkable that the dynamics of the 
martensitic transition is still to date an outstanding open modelling problem \cite{JB02}, \cite{Muller}.
%

The modelling approach we propose for the dynamics of the AM transformation consists of a fragmentation 
process in which the microstructure evolves via formation of thin plates of martensite embedded in a 
medium representing the austenite. Here, the plates are idealized as planar surfaces which nucleate and 
propagate according to a stochastic process in space and time and that are allowed to grow parallel to a 
given set of admissible directions. Our main focus will be on a simple model in two dimensions and
we imagine that the unit square represents the reference configuration occupied by the material in the 
austenitic phase. In our caricature model, nucleation events occur as a Poisson process in this domain 
and as each point appears it generates, with equal probability, a plate which propagates orthogonally 
in either the $x$ or $y$ direction. 
These plates propagate until they hit a part of the structure of plates that has already appeared or the 
boundary of the domain and as result a complicated structure of rectangles soon emerges.  A physical 
process which generates patterns of the type we are interested in is shown in Figure~\ref{1609181210} 
and the initial steps of a realization of the model are shown in Figure~\ref{cartoon}.

It is believed that the acoustic emissions in experiments on the AM transformation 
are associated with the formation of plates and the power laws observed should correspond to power laws
in the number of plates. By encoding our probabilistic model into a general branching random walk, 
we are able to obtain predictions for the asymptotic behaviour of features of the random pattern 
generated. The question that we address is the behaviour of the number of plates of a given size. We 
treat two versions of this problem. 
\begin{enumerate}
\item We look at the asymptotics for the sizes of small interfaces and show that this depends on
the way the pattern is generated and is not a power law.
\item For a complete fragmentation we consider the number
of interfaces above a certain size and obtain a power law for this quantity.
\end{enumerate}

An analogous modelling approach for martensitic microstructure has been adopted
for the first time, to the best of our knowledge, in \cite{Pasko}, where the structural properties of 
martensite are deduced from the Fourier analysis of a 2D pattern obtained from a fragmentation process.
We also refer to \cite{Torrens16} for an alternative probabilistic model based on a fragmentation scheme
which takes place over a finite-spacing lattice and addresses our second question.  Besides the 
discrete nature of the model \cite{Torrens16}, 
the main difference with the present investigation lies in the mathematical approaches adopted for the 
analysis. In \cite{Torrens16} the investigation is focused on the exact solution of a complicated 
finite-difference equation obtained `from scratch' for the conditional probability of splitting 
events occurring at the lattice level and the outcome is limited to the construction of the asymptotic 
mean profile of the geometric quantities of interest. Conversely, in the present scenario general 
branching process theory provides finer information on the shape of the asymptotic distribution of 
the interfaces as it shows there is almost sure convergence of the quantities of interest to a limit 
random variable. It also allows us to vary the model more easily.

There are many papers which discuss models for patterns created by fragmentation. This model is referred
to as the Gilbert tesselation in \cite{MM}, where some features of the patterns are discussed.
For mathematical work on 
the analysis of other fragmentation processes of relevance here, see  \cite{Ber} and  \cite{HKK}. In 
the scalar case, we refer to~\cite{Derrida87} and~\cite{Frontera95} where a simple fragmentation model 
has been introduced 
for the description of first order phase transitions occurring as avalanches. For two dimensional 
fragmentation processes there is work on quadtrees which provides a more sophisticated analysis \cite{BNS} 
of a model related to the one we discuss here. Alternative approaches may be found in work on random 
tessellations. For tessellations of a compact region generated by iterations we refer to \cite{NW}, and 
branching wise \cite{GST}.

%
%
%
%


\subsection{The model}\label{sec:model}

\begin{figure}[h!]
    \centering
    \begin{subfigure}[b]{0.2\textwidth}
\includegraphics[width=1.8cm]{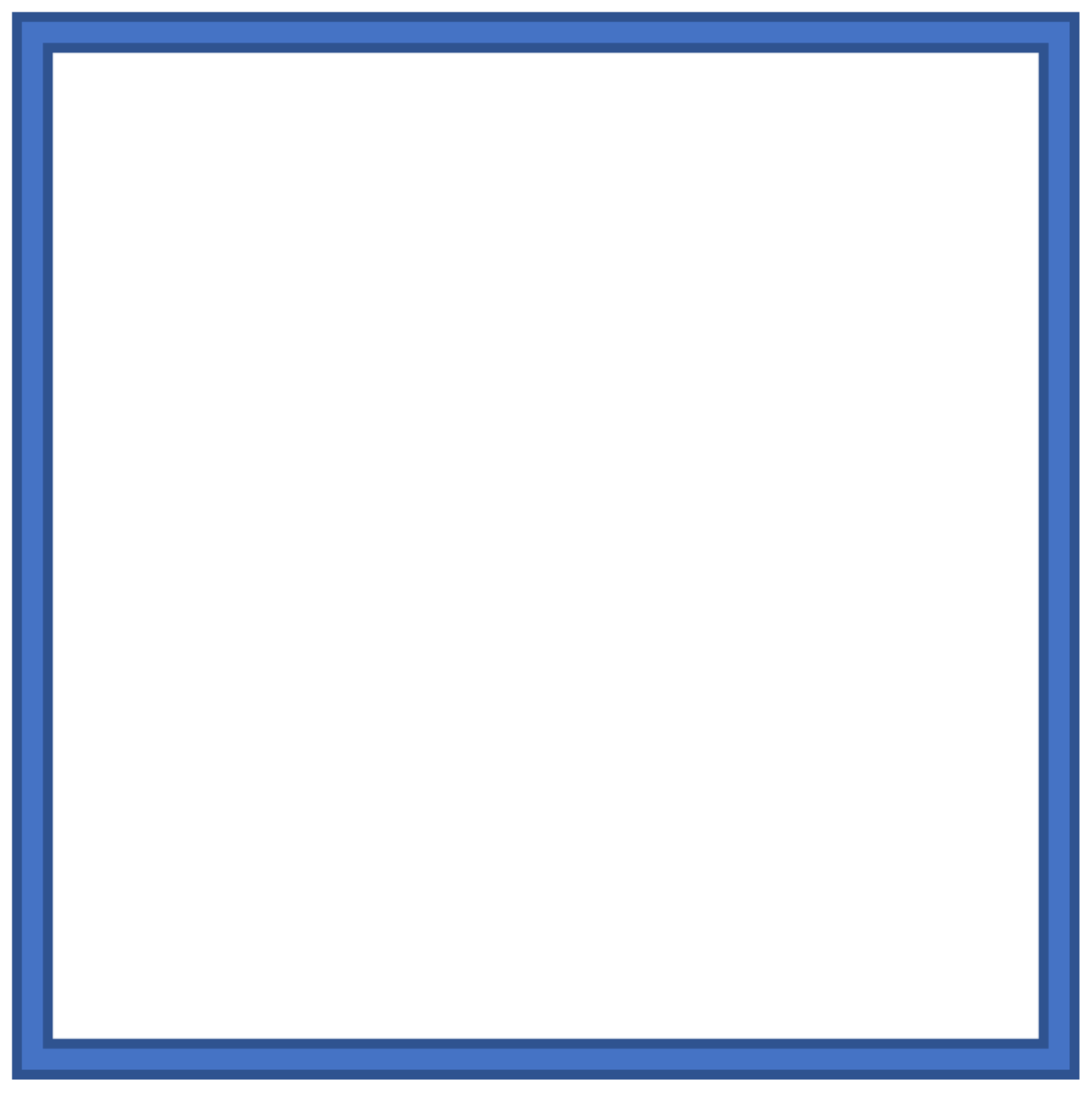}
        \caption{}\label{cartoon_a}
    \end{subfigure}
\hspace{-0.9cm}
    ~ 
    \begin{subfigure}[b]{0.2\textwidth}
\includegraphics[width=1.8cm]{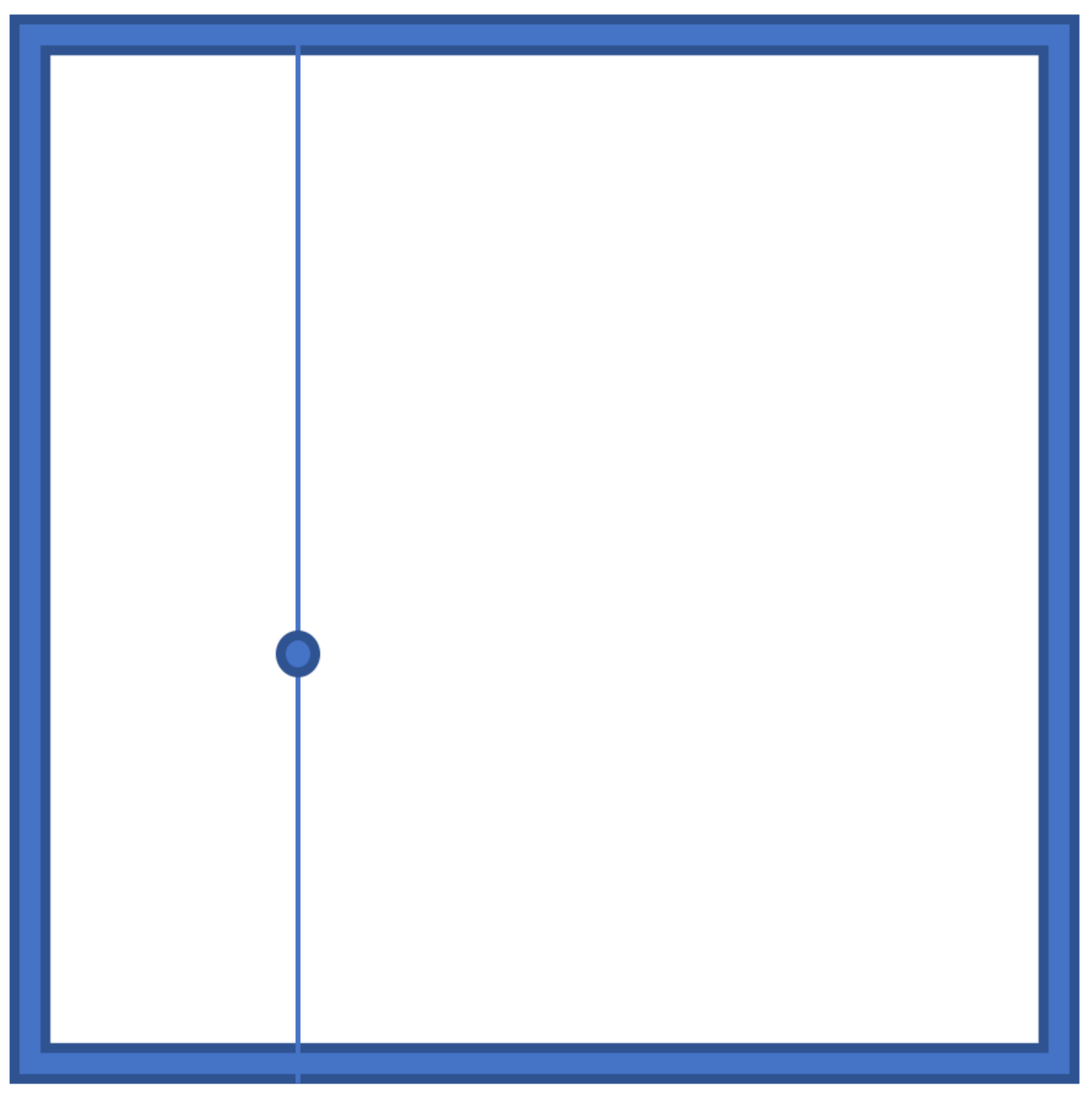}
        \caption{}\label{cartoon_b}
    \end{subfigure}
\hspace{-0.8cm}
    ~ 
    \begin{subfigure}[b]{0.2\textwidth}
\includegraphics[width=1.8cm]{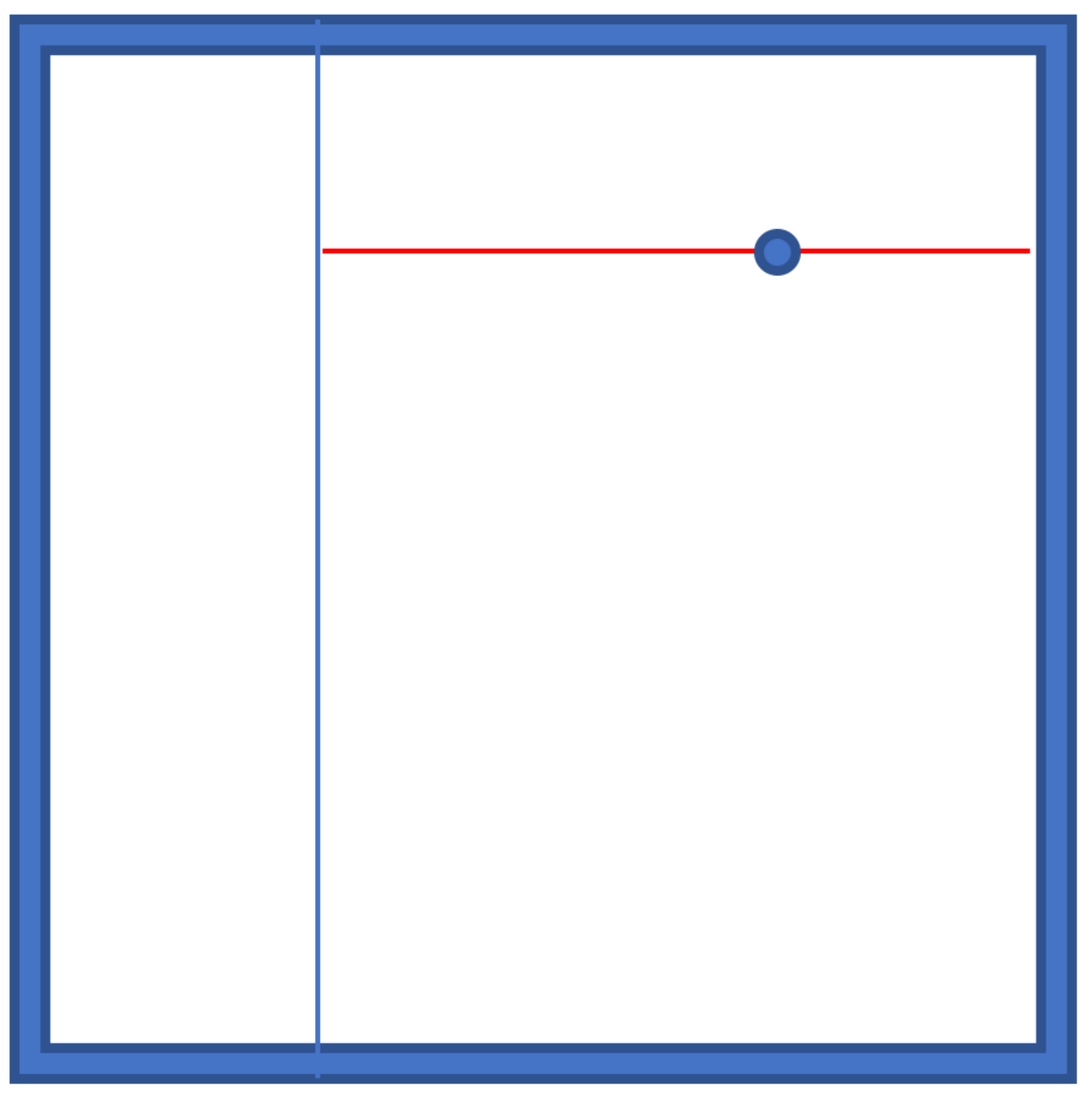}
        \caption{}
        \label{cartoon_c}
    \end{subfigure}
\hspace{-0.7cm}
    \begin{subfigure}[b]{0.2\textwidth}
\includegraphics[width=1.8cm]{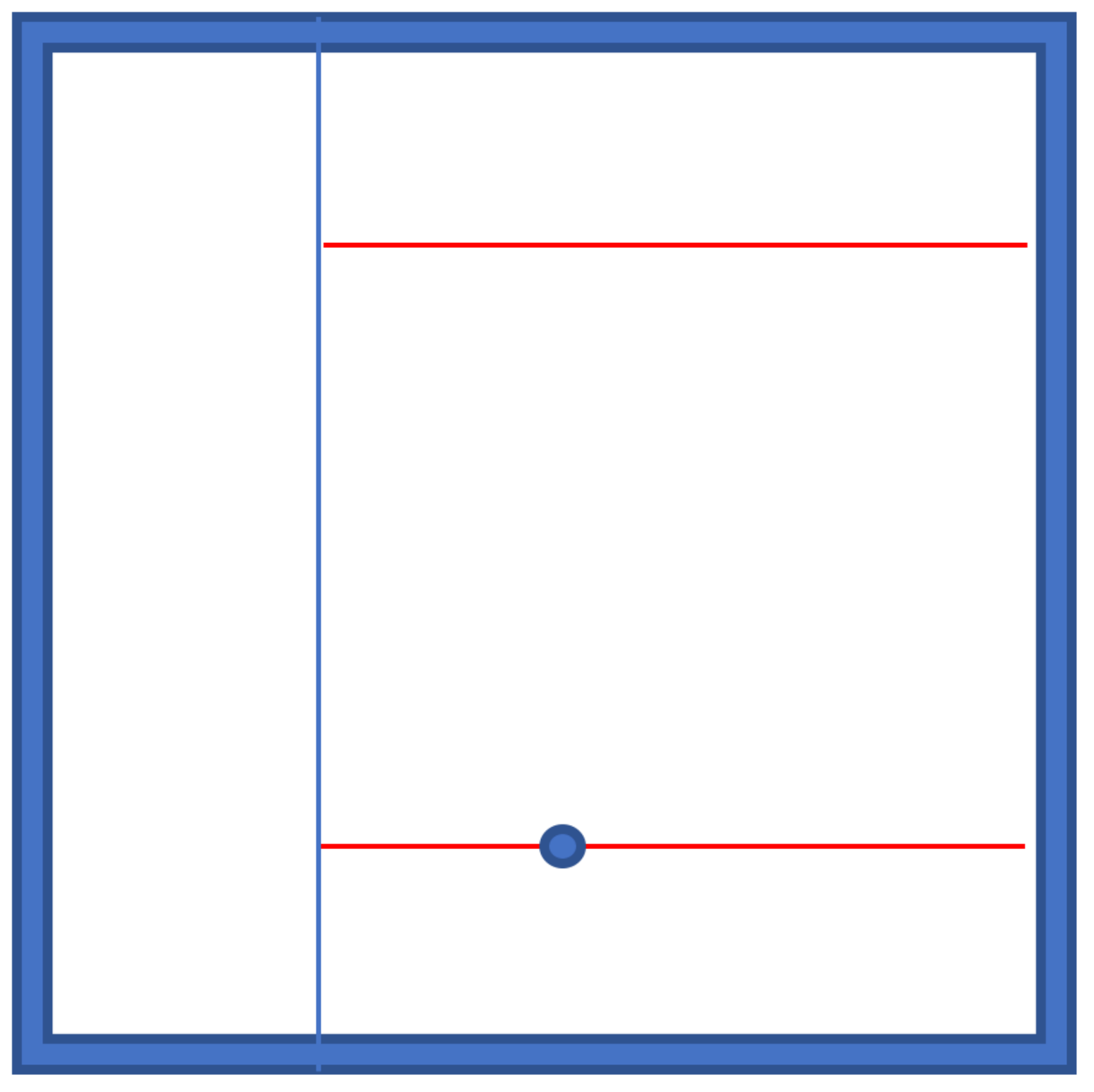}
        \caption{}
        \label{cartoon_d}
    \end{subfigure}
\hspace{-0.7cm}
    \begin{subfigure}[b]{0.2\textwidth}
\includegraphics[width=1.8cm]{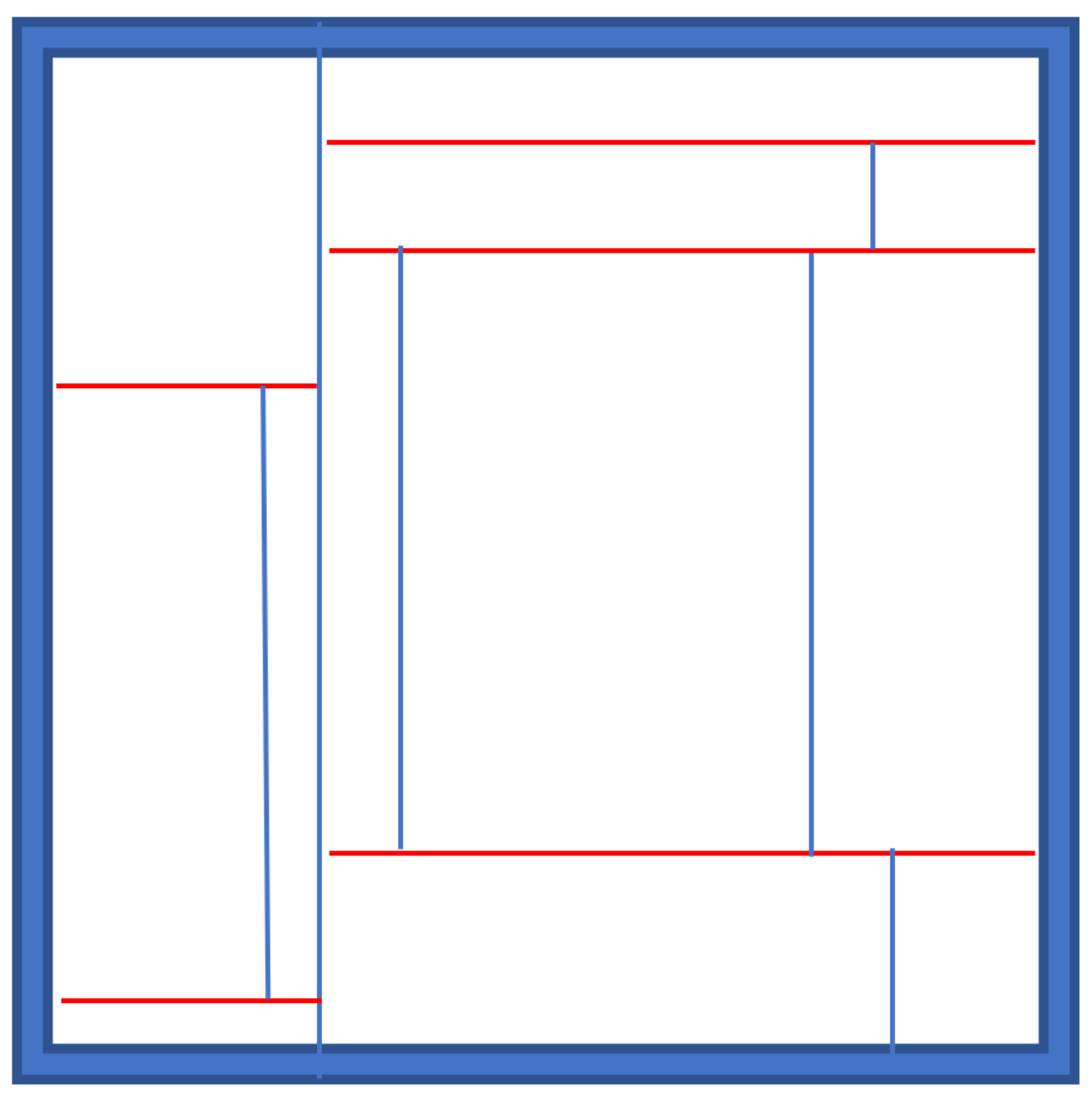}
        \caption{}
        \label{cartoon_e}
    \end{subfigure}
    \caption{(a)-(d) show the initial stages of a fragmentation process. (e) is the configuration after 11 
    nucleation events
    }\label{cartoon}
\end{figure}

The model divides the unit square up into a succession of rectangles which then evolve in time to get broken up 
into smaller and smaller pieces. The initial reference configuration is a unit square (Figure~\ref{cartoon_a})
representing a homogeneous lattice in the high-symmetry phase (austenite). A point is selected in the square 
according to a uniform probability distribution and a line through the point is drawn in one direction, 
either horizontal (with probability $p$) or vertical (with probability $1-p$). The line extends through the 
square until it hits the boundary of the domain and as a result the initial square is fragmented into two
rectangles (Figure~\ref{cartoon_b}). Now we iterate the procedure by picking a second point, again chosen 
according to a uniform probability distribution in the initial unit square and draw a line passing through 
this point in either direction, horizontal or vertical with probability $p$ and $1-p$ respectively (in the
illustration in Figure~\ref{cartoon_c}, the horizontal direction). Now the line extends through the domain 
until it hits the boundary of the unit square or a part of the domain which has already transformed.
 
 %
 %

We continue iterating this procedure by choosing a nucleation point at random in the unit square and splitting the rectangle containing it, soon resulting in a complicated pattern (Figure~\ref{cartoon_d} and~\ref{cartoon_e}). The
requirement that interfaces cannot overlap is an important ingredient of the model inspired by experiments. 
It constitutes the modelling mechanism for a material where the parts which have already transformed into 
martensite cannot evolve any further.
\begin{figure}[h!]
\centering%
\includegraphics[width=5.5cm]{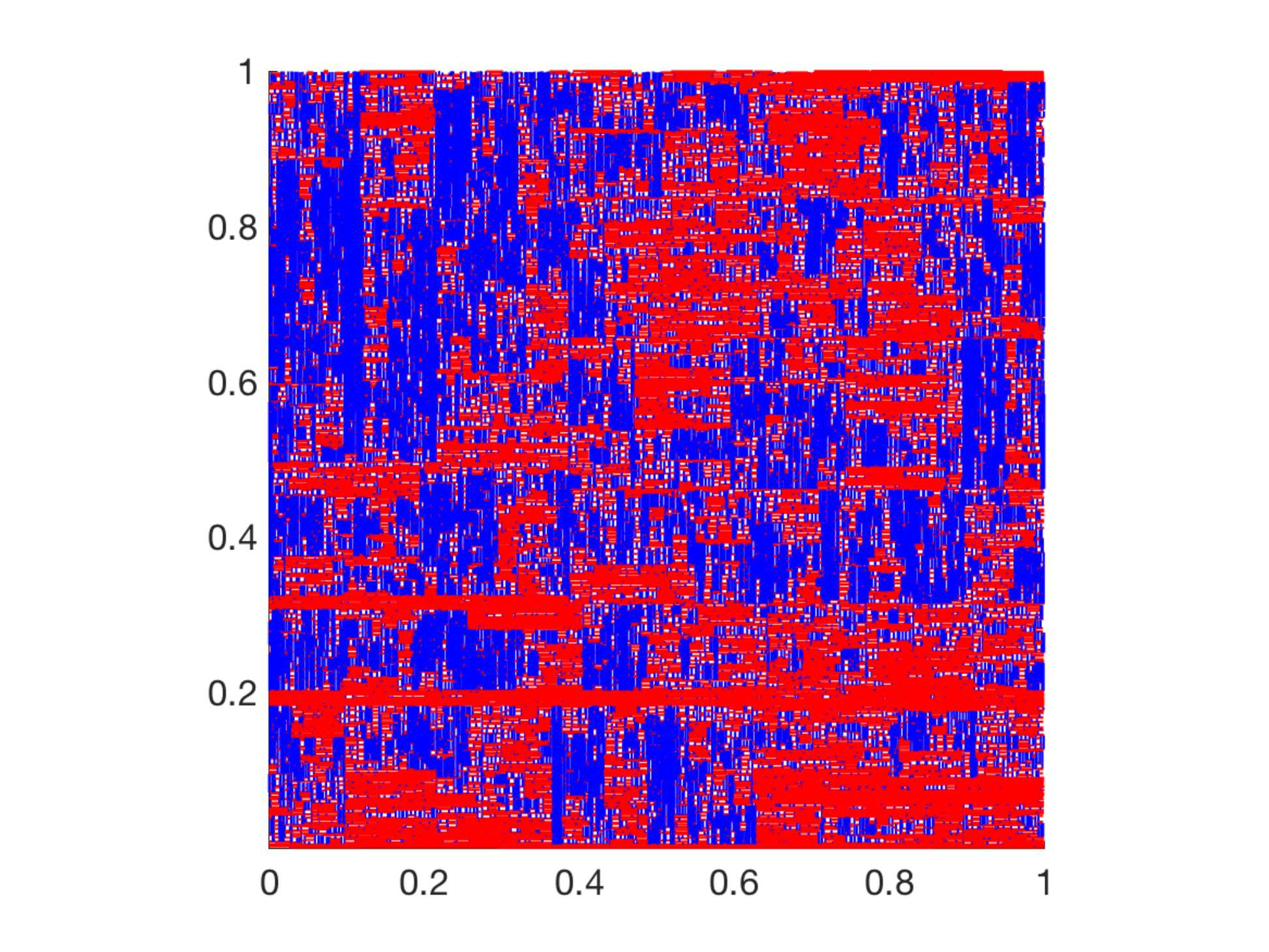}
\hspace{-1.1cm}
\includegraphics[width=5.5cm]{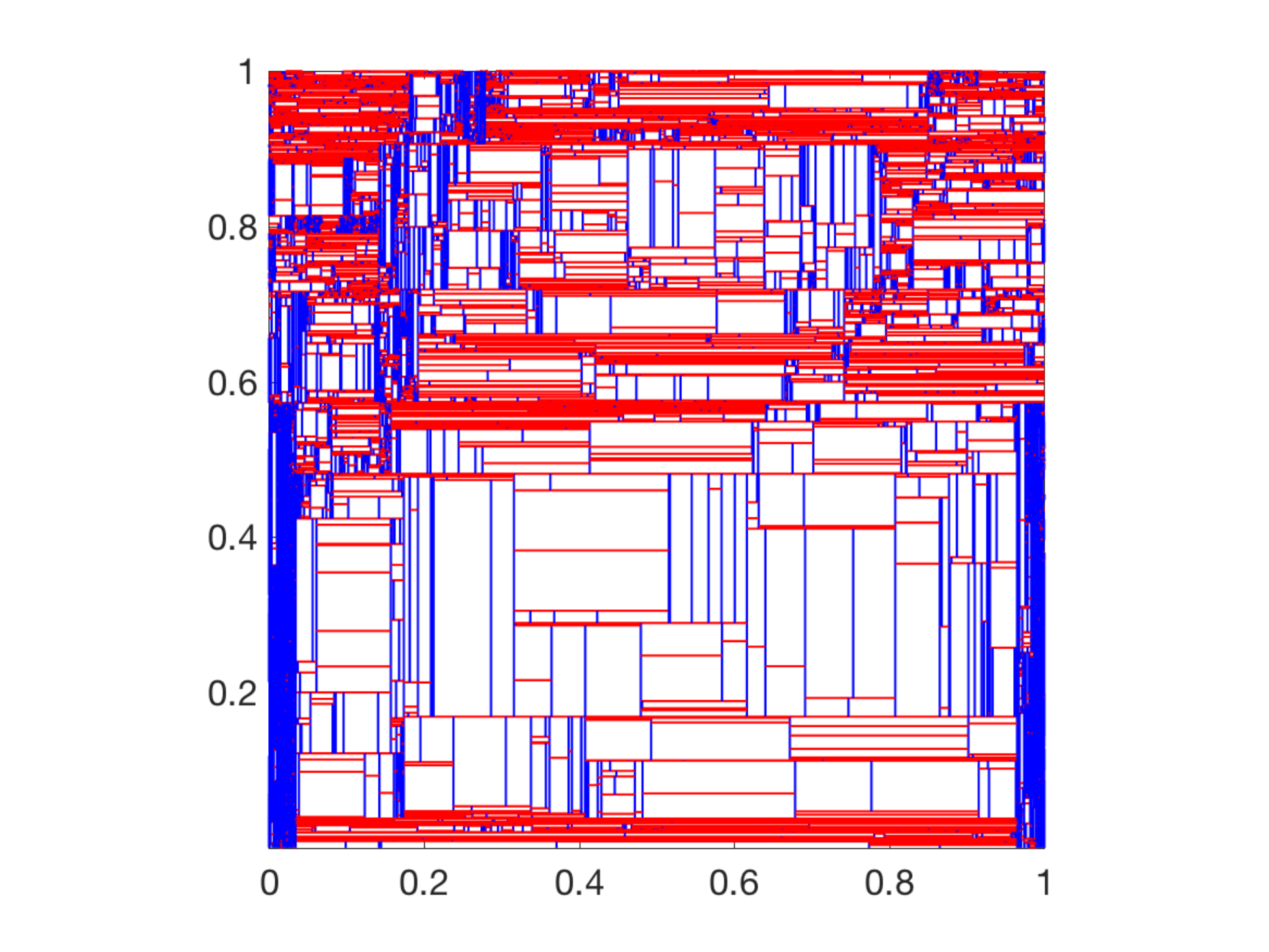}
\caption{
Left: microstructure obtained for the fragmentation model of Section \ref{sec:model}. At each step we 
split the rectangle of largest area. Here $p=\frac{1}{2}$.
Right. Pattern obtained when each existing rectangle is split at the same time.
Both microstructures are composed of $2^{15}$ rectangles.
}\label{2006172133}
\end{figure}
The procedure described above is equivalent to requiring the nucleation to occur in a rectangle with 
a probability proportional to its area. This process can be continued indefinitely to get a complete 
fragmentation of the unit square. It is easy to see that, as each rectangle is split in two and the 
future evolution inside a rectangle is independent of what occurs outside, the individual rectangles 
that arise can be encoded by an infinite binary tree and then any pattern we see at a finite scale can 
be obtained by taking a cut through the branches of this tree. We will be interested in different ways 
of looking at this tree and will consider two time parametrizations that lead to a tractable analysis. 
Simulations of the outcome from the two different versions are shown in Figure~\ref{2006172133}.


In Figure~\ref{2006172133}-Left we show a realization of the pattern obtained when, at each 
iteration, a new nucleation point is selected according to a uniform probability distribution in 
the rectangle of largest area. 
In order to do this we introduce a time parametrization so that at time $t$ in this evolution 
all rectangles in the pattern 
have area at most $e^{-t}$. This is the fragmentation 
mechanism analysed in Section \ref{sec:largest}.

Figure~\ref{2006172133}-Right shows one realization of the fragmentation process described in 
Section~\ref{sec:indep} which consists of splitting, at a fixed time, all the rectangles available 
regardless of their size.  This fragmentation mechanism
leaves $2^n$ rectangles after $n$ steps,
resulting in a pattern which after a finite number of steps is highly inhomogeneous. 
This is evident by observing the microstructure of Figure~\ref{2006172133}-Left which is composed of the 
same number of rectangles as Figure~\ref{2006172133}-Right.
This difference is made rigorous in the theoretical analysis in Section~\ref{1801071839} and empirically 
through the statistics of the distribution of rectangles in Section~\ref{1708281611}.

In order to describe the model mathematically we consider each rectangle as determined by two coordinates 
$(a,b)$ which are the side lengths for the edges parallel to the horizontal and vertical axes respectively. 
Nucleation events occur as a Poisson process in time and lead to the formation of an interface. 
The nucleation point is uniformly distributed in the rectangle and with probability $p$ the interface
grows horizontally or with probability $1-p$ it grows vertically
 until it reaches the edges of the rectangle in 
the vertical or horizontal direction respectively. In this way we can see that the rectangle $(a,b)$ evolves as
\[ (a,b) \to \left\{\begin{array}{ll} \{(a,bU), (a,b(1-U)) \} & \mbox{ probability } p \\  \{(aU,b), (a(1-U),b) \} 
& \mbox{ probability } 1-p \end{array} \right. \]
where $U$ is a uniform random variable on $[0,1]$.

%
%

\begin{figure}[h!]
\centering%
\includegraphics[width=5.7cm]{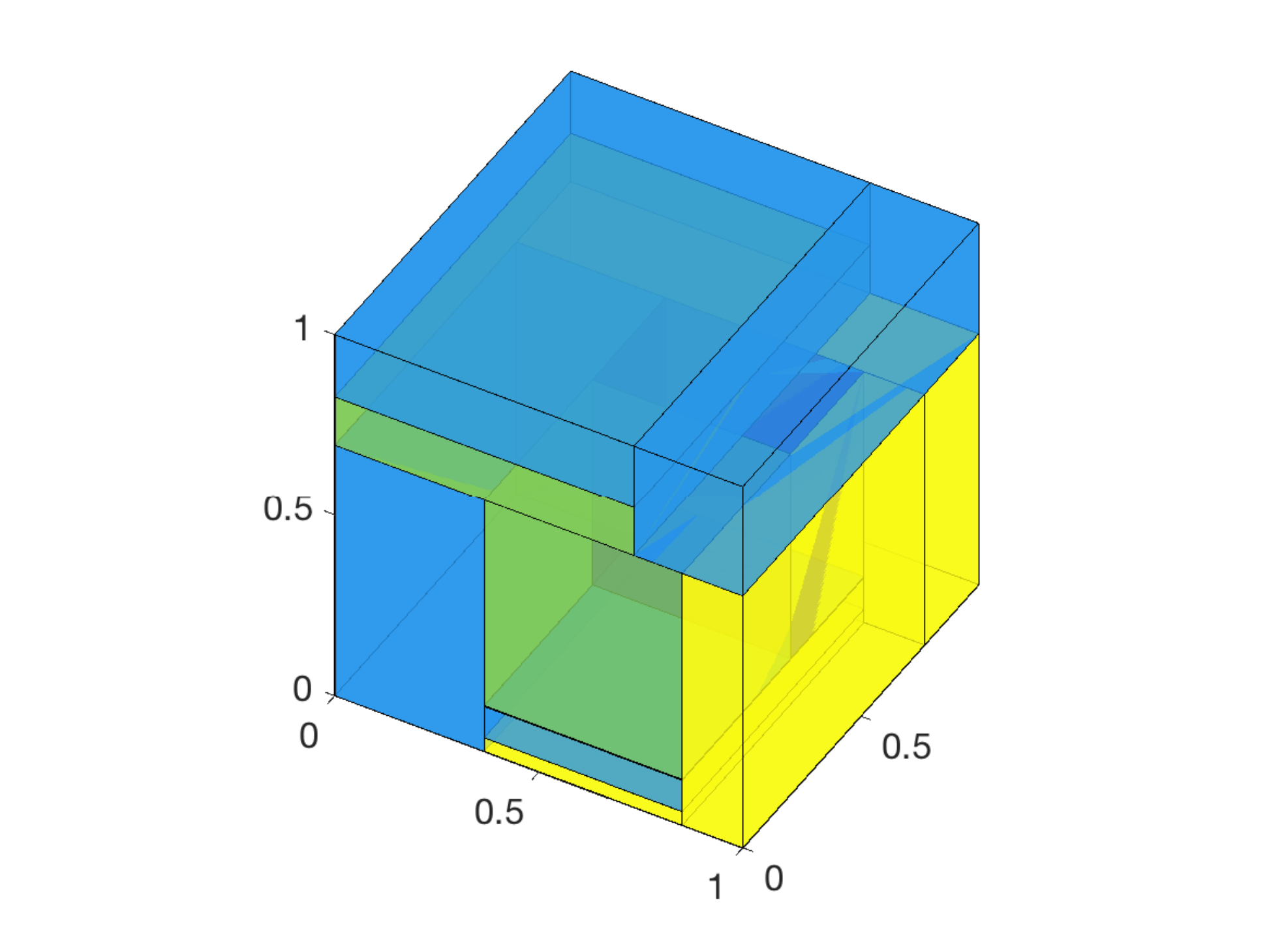}
\hspace{-1.2cm}
\includegraphics[width=5.7cm]{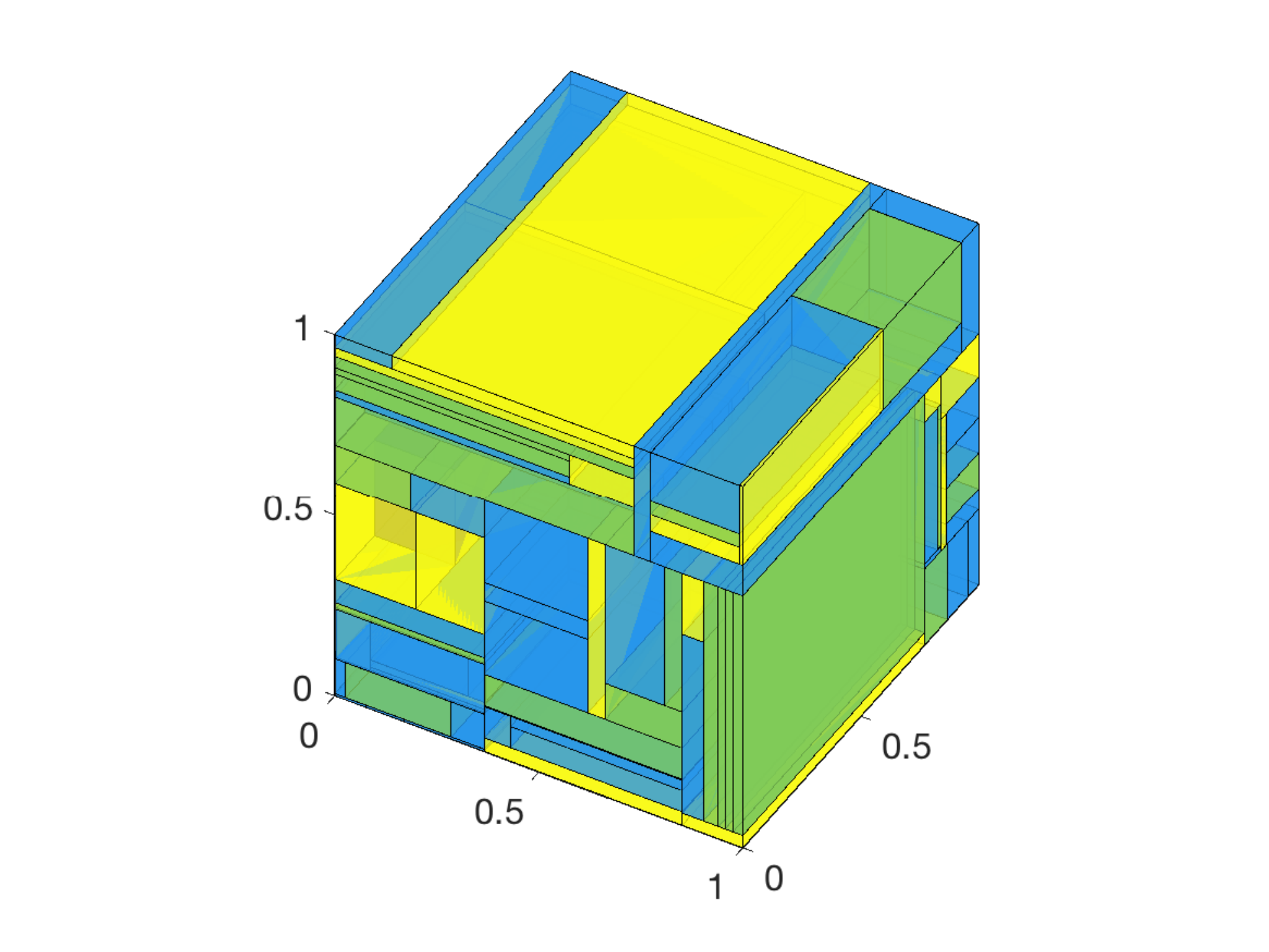}
\hspace{-1.2cm}
\includegraphics[width=5.7cm]{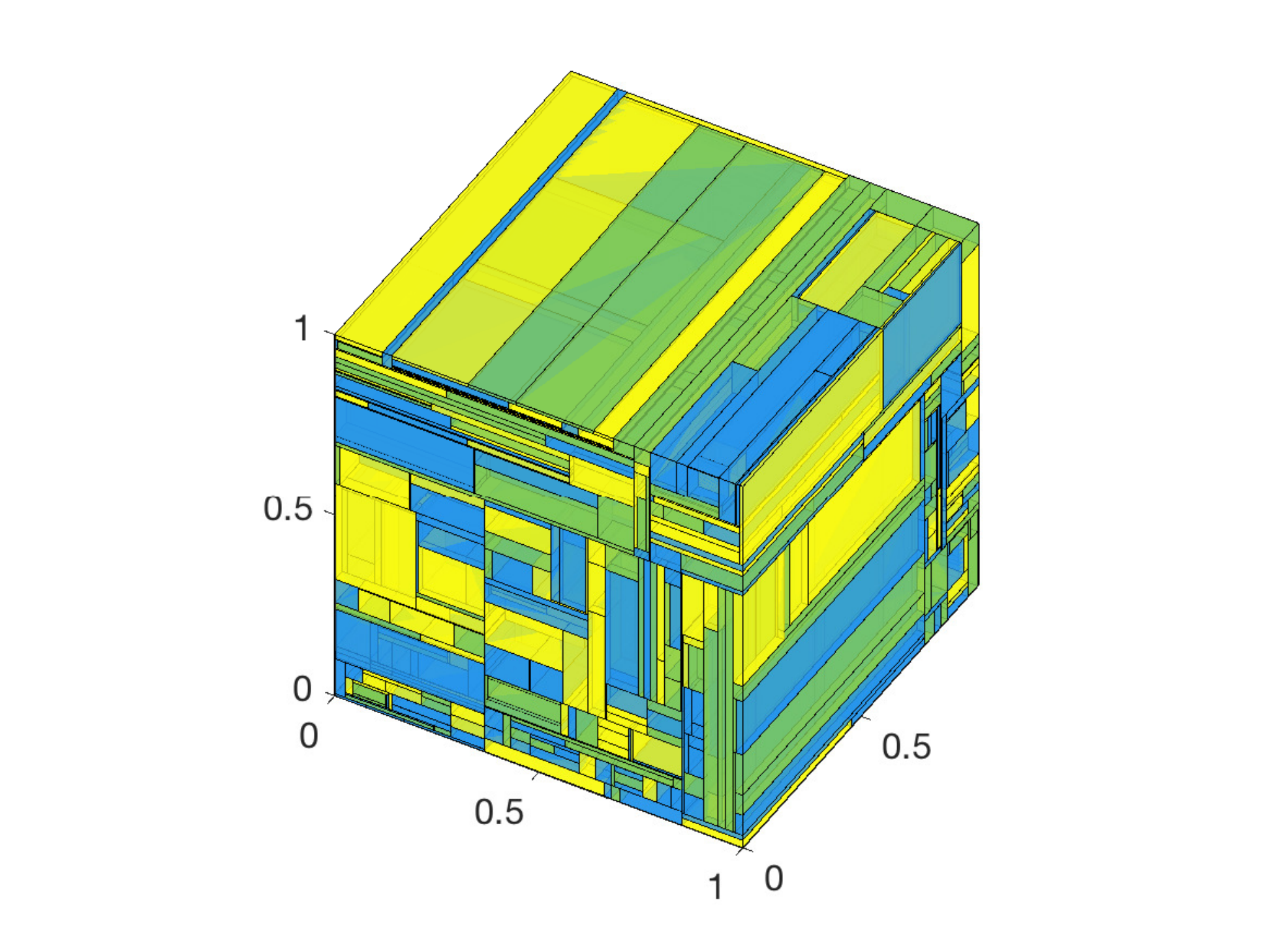}
\caption{
Three stages of the 3D fragmentation model observed after $n$ iterations. From 
left to right, $n=10,100,1000$). Here $p_1=p_2=\frac{1}{3} $. 
}
\label{cubes}
\end{figure}

The three-dimensional version of the fragmentation scheme is obtained analogously to the two-dimensional 
version by breaking down a unit cube
into many rectangular cuboids by drawing planar plates.
%
%
The procedure starts when a point is picked according to a uniform probability distribution on the unit 
cube. A plate through the point and orthogonal to either the $e_1,e_2$ or $e_3$ direction (assuming the 
standard Euclidian basis) with probability respectively $p_1,p_2$ or $1-p_1-p_2$ (with $0\leq p_1,p_2,
p_1+p_2\leq 1$) splits the cube into two rectangular cuboids. 
We iterate this procedure by picking another point in the interior of the rectangular cuboid of largest 
volume and repeating the fragmentation procedure by allowing the plates to grow till they 
hit the boundary of the cube or another part of the domain which has already transformed 
(see Figure~\ref{cubes}).
 
\subsection{The main results}

Our aim is to look at certain asymptotic quantities associated with the resulting collection of rectangles. 
In particular we would like to determine the power law growth of the number of interfaces of at least a 
certain size. We now discuss two ways of looking at the model and give two theorems that show the different 
type of asymptotic analysis that can be performed.

Firstly we consider the asymptotics of the sizes of small rectangles in the pattern when we split the
largest area first.
Let $\tilde{N}_t(A)$ denote the number of rectangles with side length $(a,b)\in A\subset [0,1]^2$ at time $t$. 
We will show the following:
\begin{thm}\label{1801071705}
Let $S$ be a closed convex set with non-empty interior in $[0,1]^2$. 
\begin{enumerate}
\item If $S\cap\{(a,b): ab=e^{-1}\}=\emptyset$, then
\[ t^{-1} \log{\tilde{N}_t(e^{-t}S)} \to -\infty, \;\; t\to\infty, \;\; a.s. \]
\item If $S \cap\{(x,y):ab=e^{-1}\}\neq \emptyset$, then setting $\beta = \sup\{\gamma^*(a,b):(a,b)\in S\}$, where 
\[ \gamma^*(a,b) = 2\sqrt{\log{a}\log{b}}, \;\; ab=e^{-1}, \]
and $\gamma^*(a,b)=0$ for all other values of $a,b$, we have 
\[ t^{-1}\log{\tilde{N}}_t(e^{-t}S) \to \beta, \;\; t\to\infty\;\; a.s. \]
\end{enumerate}
\end{thm}
In the case where each rectangle splits independently at the same constant rate we obtain a different 
function $\gamma^*$ which is strictly positive over an open set and given in \eqref{eq:gamma-indep}.

For our second result we do not consider time but look at the number of interfaces that have ever appeared 
with length greater than say $x$. In the two-dimensional version of the model this quantity is only finite 
if there is a bias in the direction of interface formation. However, in the three dimensional unbiased 
version of the model, in which we split rectangular cuboids uniformly in the same way as above, we obtain 
a power law exponent. More precisely, if we set $N_h(x)$ to be the number of horizontal plates with area 
at least $x$ in the limit pattern, then 
\begin{thm}\label{thm:3dP}
There exists a non-degenerate random variable $W>0$ such that
\[ \lim_{x\to 0} N_h(x)x^3 = W, \;\; a.s. \]
\end{thm}
By symmetry the same result holds for the plates which are parallel to the other axial directions and hence the total number of plates of more than a given area is power law with exponent 3. 
As this is a limit result, the initial domain is not important and we would expect to see that any 
material which fragmented uniformly at random by plates orthogonal to the axial directions would have 
the exponent 3 for decay of the areas of the interfaces.

The structure of the paper is as follows. We begin by describing the general branching process and the 
general branching random walk along with the main theorems concerning the growth rates of these processes 
in Section~2. In Section~3 we show how our model can be encoded as a general branching random walk. Using 
this we perform an analysis of the asymptotics for the number of rectangles of a given size within the 
pattern as a proxy for the lengths of the interfacial plates. In the case where we split the largest area 
first, this requires an extension of previous results on the general branching random walk. We then 
consider an alternative formulation in Section~4 in which we can calculate the asymptotics for the number 
of edges greater than a given length but only in a biased version of the model. In Section~5 we consider 
a variety of extensions of the model, including the three-dimensional case, and illustrate some further analytic results that can be obtained. 
In Section~6 we perform numerical experiments which confirm the analytical results and discuss the issues 
which need to be addressed when recovering these asymptotic results via simulations.

%
%
%
%
 
\section{General branching processes}\label{1801071836}

In the next section we will code up our fragmentation problem using a general branching random walk. 
The whole structure can be captured in a labelled tree as the evolution inside any rectangle does not 
affect what happens outside that rectangle. Thus in our tree each vertex will represent a rectangle 
and we put labels on the vertices to capture the information that we need about the size of the 
rectangle. The natural setting for this is the theory of general branching processes.

In the Crump-Mode-Jagers or general branching process (GBP), we model the evolution of a population 
of individuals who evolve independently. The typical individual $x$ is born at time $\sigma_x$, has 
offspring whose birth times are determined by a point process $\xi_x$ on $(0, \infty)$, a lifetime modelled 
as a non-negative random variable $L_x$, and a (possibly random) c\`adl\`ag function $\phi_x$ on $\R$ 
called a characteristic. For a general branching random walk (GBRW) we include a process for the birth 
position $\eta_x$ as well as the birth time. We will develop the set up for the GBRW and then, by ignoring 
the spatial structure, we will have a GBP.

We now establish the model formally.
Let $I_n = \cup_{k=0}^n \N^k$ and $I = \cup_k I_k$.
A sequence in $I$ will be denoted $\bfi \in I$, though we will write $i\in I_1$.
$I_0$ is the empty sequence $\emptyset$. We write $\bfi\bfj$ for the concatenation of sequences $\bfi$ and 
$\bfj \in I$. For a point $\bfi \in I\backslash I_n$, denote by $\bfi|_n$ the sequence truncated to length $n$
and by $\bfi[n]$ the $n$-th element of $\bfi$.
We write $\bfj \leq \bfi$ if $\bfi = \bfj \bfk$ for some $\bfk$, and denote by $|\bfi|$ the length of the 
sequence $\bfi$.

If $\bfi$ has $j$ children, then these are $\bfi 1, \bfi 2, \ldots, \bfi j$.
A tree $T$ is a subset of the space $I$ such that:
\begin{itemize}
\item[i)] $\emptyset \in T$ the root of the tree;
\item[ii)] $\bfi \in T$ implies $\bfi|_k \in T$ for all $k < |\bfi|$;
\item[iii)] if $\bfi j \in T$ for some $j \in \N$, then $\bfi 1, \bfi 2, \ldots, \bfi (j-1) \in T$.
\end{itemize}
We have a one-one correspondence between trees and feasible realisations of the set of individuals in a branching 
process. In this paper, unless otherwise stated, from the next section the tree will be fixed as a binary tree.

We now introduce suitable labels on the tree. 
For each element $\bfi \in T$, let the life-story for $\bfi$ be
$V_\bfi = (L_\bfi, \xi_\bfi, \eta_\bfi, \chi_\bfi)$.  
The $V_\bfi$ are i.i.d., $L_\bfi \in \R_+$ is the life span, 
$\xi_\bfi : \R_+ \to \Z_+$ (non-decreasing) is the birth time process, 
$\eta_\bfi : \R_+ \to \R^d$ is the birth position process, and
$\chi_\bfi : \R_+ \to \R_+$ is the characteristic, a function used to count attributes of the individual.
For $t \leq L_\bfi$, $\xi_\bfi(t)$ is the number of children born to $\bfi$ up to and including time $t$.  
Let $N_\bfi = \xi_\bfi(L_\bfi)$ be the total number of children born to $\bfi$.  
For a general branching random walk, we make no assumption about the joint
distribution of $L_\bfi, \xi_\bfi, \eta_\bfi$ and $\chi_\bfi$, though we can w.l.o.g.\ assume 
that $\xi_\bfi(t) = \xi_\bfi(t \wedge L_\bfi)$.  
Let the birth time of $\emptyset$ be
$\si_\emptyset = 0$, and let the birth time for $\bfi j$ be $\si_{\bfi  j} = \si_\bfi + t_\bfi(j)$, 
where $t_\bfi(j) = \inf \{ t \,:\, \xi_\bfi(t) \geq j \}$.  Any distribution for $V_\emptyset$ induces a
probability space $(\Om, \calB, \Pb)$, where $\Om$ is the space of
trees $T$ and associated life-stories $\{ V_\bfi \,:\, \bfi \in T \}$.

\subsection{The general branching process}

A key result for the GBP was obtained by Nerman \cite{Ner81} in
which he established a reasonably general strong law of large numbers for the GBP.

For the GBP we take $\eta_\bfi=0$ in our earlier framework, so that we ignore the spatial component of the 
process. For the general branching process individual $\bfi$ has $\xi_\bfi(0, \infty)$ offspring whose birth 
times $\sigma_{\bfi i}$ satisfy
$$
\xi_\bfi = \sum_{i = 1}^{\xi_\bfi(\infty)} \delta_{\sigma_{\bfi i} - \sigma_\bfi}.
$$
The trace of the underlying Galton-Watson branching process is a random subtree of $I$ which we denote 
by $\Sigma$.

We assume that the triples $(\xi_\bfi, L_\bfi, \chi_\bfi)_\bfi$ are i.i.d. but make no assumptions 
on the joint distribution of $(\xi_\bfi, L_\bfi, \chi_\bfi)$. We write $(\xi, L, \chi)$ for an unspecified
individual. We will write $\bp$ for the law of the general branching process and $\be$ for its expectation.

Let
\[ \xi(t) = \xi((0, t]), \quad \nu(dt) = \be \xi(dt), \quad \xi_\gamma(dt) = e^{-\gamma t} \xi(dt), 
\quad \text{and} \quad \nu_\gamma(dt) = \be \xi_\gamma(dt), \]
for $\gamma \in (0,\infty)$. We assume that the general branching process is super-critical in that there 
exists a Malthusian parameter $\alpha \in (0, \infty)$ given by the value of $\gamma$ such 
that $\nu_\gamma(\infty) = 1$.  

We will also need $\mu_1 = \int_0^{\infty} t \nu_{\alpha}(dt)$, the first moment of the probability 
measure $\nu_\alpha$.

The characteristic $\chi$ enables us to count quantities in the population. The characteristic counting process
$Z^\chi$ is defined as
\begin{equation}
Z^\chi(t) = \sum_{\bfi \in \Sigma} \chi_\bfi (t - \sigma_\bfi).
\end{equation}
There is a recursive decomposition of $Z^\chi$
\begin{equation}
Z^\chi(t) = \sum_{\bfi \in \Sigma} \chi_\bfi (t - \sigma_\bfi) = \chi_\emptyset(t) + 
\sum_{i = 1}^{\xi_\emptyset(\infty)} Z_i^\chi(t- \sigma_i), \label{eq:ccp}
\end{equation}
where the $Z^\chi_i$ are i.i.d.\ copies of $Z^\chi$. 

To state the strong law of large numbers for the process $Z^{\chi}$ we require two more quantities. 
Firstly we define the discounted mean process and the discounted characteristic
$$
z^\chi(t) = e^{-\alpha t} \be Z^\chi(t) \quad \text{and} \quad u^\chi(t) = e^{-\alpha t} \be \chi(t).
$$
It is easy to check that these satisfy the renewal equation
\begin{equation}
z^\chi(t) = u^\chi(t) + \int_0^\infty z^\chi(t-s) \nu_\alpha(ds). \label{eq:zrenewal}
\end{equation}

The second quantity we require is the analogue of the classical branching process martingale defined by
\[ M_t = \sum_{\bfi \in \Lambda_t} e^{-\alpha \sigma_\bfi}, \]
where
\[ \Lambda_t = \{\bfi\in\Sigma:\:\bfi = \bfj i \mbox{ for some }\bfj\in\Sigma,i\in\bn,\mbox{ and }\sigma_{\bfj}
\leq t < \sigma_\bfi\} \]
is the set of individuals born after time $t$ to parents born up to time $t$. The process $M$ is a 
non-negative c\`adl\`ag $\cf_t$-martingale with unit expectation, where
$$
\cf_t = \sigma(\cf_\bfi, \sigma_\bfi \leq t) \quad\text{and} \quad \cf_\bfi = \sigma(\{(\xi_\bfj, L_\bfj): 
\sigma_\bfj \leq \sigma_\bfi\}).
$$

As $M$ is a positive martingale, by the martingale convergence theorem, $M_t \to M_\infty$ as 
$t \to \infty$, almost-surely, for some random variable $M_\infty$. Furthermore, under the condition that 
$\be \left[\xi_\gamma(\infty) (\log \xi_\gamma(\infty))_+ \right] < \infty$, the martingale $M$ 
is uniformly integrable, there is $L^1$ convergence $M_t\to M_{\infty}$ and the limit random variable 
is non-degenerate in that $0<M_{\infty}<\infty$.

We now give Nerman's Theorem, a strong law for the general branching process counted with characteristic $\chi$. This version is Theorem~2.5 in \cite{CCH} (coupled with the comments in the paragraph 
preceding the statement of the Theorem).

\begin{thm}\label{thm:sllngbp}
Let $(\xi_\bfi, L_\bfi, \chi_\bfi)_\bfi$ be a general branching process with Malthusian parameter $\alpha$, 
where $\chi \geq 0$ and $\chi(t) = 0$ for $t < 0$. Assume that $\nu_\alpha$ is non-lattice. 
Assume further that $\be \xi(\infty)<\infty$ and there exists a non-increasing bounded positive 
integrable c\`adl\`ag function $h$ on $[0, \infty)$ such that
$\be \left( \sup_{t \geq 0} \frac{e^{-\gamma t} \chi(t)}{h(t)} \right) < \infty$.
Then,	
$$ z^\chi(t) \to z^\chi(\infty) = \mu_1^{-1}{\int_0^\infty u^\chi(s) ds}, $$
and
$$ e^{-\gamma t} Z^\chi(t) \to z^\chi(\infty) M_\infty, \text{ a.s.}, $$
as $t \to \infty$, where $M_\infty$ is the almost sure positive limit of the fundamental martingale of the 
general branching process. Furthermore, if $\be \left[\xi_\gamma(\infty) (\log \xi_\gamma(\infty))_+ \right] 
< \infty$, the convergence also takes place in $L^1$.
\end{thm}

\subsection{General branching random walks}\label{sec:gbrw}


To discuss the GBRW we introduce a spatial component to the process.
Let the birth position of $\emptyset$ be $z_\emptyset = 0$, and let the birth position of $\bfi j$ be 
$z_{\bfi j} = z_\bfi + \eta_\bfi(t_\bfi(j))$. 
Let $\de_x$ be the point mass at $x$.
It is often useful to think about $Z_t$, the point process of individuals alive at time $t$, that is
\[ Z_t =  \sum_{\bfi \in T} \chi_{\bfi}(t-\sigma_{\bfi}) \de_{z_{\bfi}}, \]
where the characteristic $\chi_{\bfi}(t) = I_{\{t<L_{\bfi}\}}$. Another useful process is
$Y_t$ the point process in $\mathbb{R}^d$ given by the coming generation at time $t$.
That is $Y_t = \sum_{\bfi j \in T, \si_\bfi < t, \si_{\bfi j} \geq t} \de_{z_{\bfi j}}$, 
Here we give a theorem about the growth and spread of such processes.

It is a classical result that a supercritical branching process will grow exponentially and that a supercritical 
branching random walk will have an associated shape theorem which indicates the region in which the 
number of particles will grow exponentially. Here we quote a result from Biggins \cite{Big94} in the $d$-dimensional case. A full version of the $1$-dimensional case is stated and proved in \cite{Big95}.

Firstly we need to define some terms. The moment generating function for the positions and birth times
of the offspring is
\[ m(\bftheta,\phi) = E \int e^{-\bftheta. \bfx-\phi t} \eta_{\emptyset}(d\bfx)\xi_{\emptyset}(dt), 
\;\; \bftheta\in \R^d, \phi \in \R_+. \]
We assume that $m(\bftheta,\phi)<\infty$ in a neighbourhood of the origin.
We set 
\[ \alpha(\bftheta) = \inf\{ \phi: m(\bftheta,\phi) \leq 1\}, \]
and $\alpha^*$ to be the Legendre transform defined by
\[ \alpha^*(\bfa) = \inf_{\bftheta}\{ \bfa.\bftheta + \alpha(\bftheta)\}. \]
We write ${\cal A} := \{ a: \alpha^*(a)>-\infty\}$.

Let $N_t$ be the number of individuals at time $t$ counted according to the characteristic $\chi$.
\[ N_t = \sum_{\bfi\in T} \delta_{z_{\bfi}} \chi_{\bfi}(t-\sigma_{\bfi}), \]
for the random measure of the positions of all the individuals counted by the characteristic $\chi$ at time $t$.
Thus if $\chi_{\bfi}(t) = I_{\{\sigma_{\bfi}<t\}}$, then $N_t(A)$ is the number of individuals that 
have been in the set $A$ by time $t$. The characteristic $\chi_{\bfi j} = I_{\{\si_\bfi < t, \si_{\bfi j} 
\geq t\}}$ gives $N_t=Z_t$, the process which records the locations of the individuals in the coming generation 
at time $t$. We will see that this captures the rectangles which arise from splitting rectangles of area 
larger then $e^{-t}$, but not yet splitting ones of smaller area.

\begin{defn}
A process is said to be \emph{well regulated} if 
\[ E \sup_t e^{-\alpha(\theta) t} \chi_{\emptyset}(t) < \infty. \]
It has \emph{exponential tails} if for any unit vector $\bf{n}$ there are $\bftheta,\phi$ 
such that $m(\bf{n}.\bftheta,\phi)<\infty$.\\
We say that the process is \emph{lattice} if the spatial location of the offspring is confined to 
lie on a discrete subgroup of $\br$.
\end{defn}

\begin{thm}\label{thm:gbrw}[\cite{Big94}~Theorem~4.1]
Suppose that the process with the characteristic $\chi$ is well regulated, non-lattice with exponential 
tails. Suppose that $A$ is a closed convex set with a non-empty interior such that $A\cap {\mathcal A} 
\neq \emptyset$. Let $\beta = \sup\{\alpha^*(a): a\in A\}$.
\begin{enumerate}
\item If $\beta<0$, then for any $\gamma>\beta$
\[ e^{-\gamma t} N_t(tA) \to 0, \;\; t\to\infty, \;\; a.s. \]
\item If $\beta>0$, then
\[ t^{-1}\log N_t(tA) \to \beta, \;\; t\to\infty\;\; a.s. \]
\end{enumerate}
\end{thm}

In the problem at hand the spatial locations and birth times of offspring will take values in $\br^{d+1}_+$ 
and hence the non-lattice condition is fulfilled. The other conditions will also follow straightforwardly 
in our application. Thus we should see exponential growth of the number of individuals at points in the 
interior of the region where $\beta>0$. Unfortunately this theorem is not directly applicable as we will see 
that, for our model where we split the largest rectangle first, the function $\alpha^*(a)$ is only finite 
on a line segment in $\mathbb{R}^2$ and thus we will establish a version of this theorem for our model directly.

We also note that there are more refined theorems in less general settings, for instance in \cite{Big90} 
and \cite{Uch} and a large subsequent literature for discrete time branching processes, Markov branching 
process, branching diffusions and fragmentation processes.

\section{Asymptotics of rectangles}\label{1801071839}

In order to capture the patterns in our model for the martensitic phase transition, we encode the model as 
a general branching random walk in $\br^2_+$. We will regard each rectangle as an individual in the branching 
process and their location will be determined by the side lengths of the corresponding rectangle. 

For each vertex $\bfi$ in the tree we have a pair $(U_{\bfi},B_{\bfi})$ where $U_{\bfi}$ is a uniformly
distributed random variable over $[0,1]$ determining the random division of the 
rectangle and $B_{\bfi}$ is a Bernoulli random variable which determines whether the interface is 
horizontal or vertical, so
\[ B_{\bfi} = \left\{ \begin{array}{ll} 1 & \mbox{horizontal} \\ 0 & \mbox{vertical.} \end{array} \right. \]
We assume that it is a Bernoulli($p$) random variable for some $0<p<1$.
Using this we can write the possible outcomes for a rectangle $(a,b)$ with index $\bfi$ from a splitting 
event as the two pairs
\[ (a,b) \to \left\{ \begin{array}{ll} ((1-B_{\bfi})U_{\bfi}a,B_{\bfi}U_{\bfi}b), \\
((1-B_{\bfi})(1-U_{\bfi})a,B_{\bfi}(1-U_{\bfi})b). \end{array} \right. \]

At this stage there is no time parameter; the tree provides a description of all the fragmentations. 
We now consider two possible time parametrizations. In the first we split the rectangles to ensure 
that at a given time all rectangles are at most a certain size. In the second every rectangle evolves
independently of the rest, the splitting time chosen according to a given distribution. We will discuss
the deterministic and exponential cases in detail.

\subsection{Splitting largest areas first}\label{sec:largest}

In our first approach we ensure the rectangles in the pattern are of roughly the same size at a given time 
by taking the birth time of the individual to be determined by their size. To incorporate the branching
structure we will make a logarithmic transformation of the coordinate system. In this way all
the individuals alive at time $t$ correspond to rectangles with areas that are at most $e^{-t}$. 

In order to convert this into a GBRW we regard each rectangle as an individual with a location given by the 
logarithm of the side lengths of the rectangle. We also introduce a time coordinate in order to ensure that 
the larger rectangles split before the smaller ones and do this using the space and time distribution of 
the offspring of individual $\bfi$:
\[  \eta_{\bfi}, \xi_{\bfi} = \left\{ \begin{array}{ll} (-(1-B_{\bfi})\log{U_{\bfi}},-B_{\bfi}\log{U_{\bfi}}) & 
-\log{U_{\bfi}} \\
 (-(1-B_{\bfi})\log(1-U_{\bfi}),-B_{\bfi}\log(1-U_{\bfi})) & -\log(1-U_{\bfi}). \end{array} \right.  \]
That is the offspring distribution of spatial positions and birth times is the point process
\begin{eqnarray*}
 & & (\eta(d\bfx),\xi(dt)) \\
 & & = \left\{ \begin{array}{ll} (\delta_{(0,-\log{x})}(d\bfx),
 \delta_{{-\log{t}}}(dt))I_{\{x=t\}} + 
 (\delta_{(0,-\log{(1-x)})}(d\bfx),\delta_{-\log{(1-t)}}(dt)) I_{\{x=t\}} & p \\
 (\delta_{(-\log{x},0)}(d\bfx),\delta_{{-\log{t}}}(dt))I_{\{x=t\}} + (\delta_{(-\log{(1-x)},0)}(d\bfx),
 \delta_{-\log{(1-t)}}(dt)) I_{\{x=t\}} & 1-p . \end{array} \right. 
 \end{eqnarray*}
Thus we can view this evolution as the individual has two offspring which are both born in a location 
displaced from their parent by a shift in either the $x$ or $y$ coordinate directions. The time of their 
births is determined by how far they are displaced. The individual dies at the time of its last birth, 
giving $L_{\bfi} = \max\{-\log{U_{\bfi}},-\log{(1-U_{\bfi})}\}$.

The associated general branching random walk, denoted by $Z_t$, is the point process of 
individuals in the population at their respective locations at time $t$. Each individual 
reproduces independently according to the point process above. In order to count all 
the individuals representing rectangles that have subdivided by time $t$ we need to
use a characteristic. This we choose to be the indicator of individuals who are born after time $t$ to 
mothers born before time $t$ - that is the coming generation, denoted by $M_t$.

We now proceed to do the concrete calculations for the interface model we have here. From this we will 
be able to compute asymptotically how many individuals are in a given region at a given time and hence 
how many rectangles there are of certain sizes. 

We recall the definitions from Section~\ref{sec:gbrw}.
Firstly we compute $m(\bftheta,\phi)$. It is straightforward to calculate this from our description of 
the individual particle evolution. We let $\bftheta = (\theta_1,\theta_2)$, 
\begin{eqnarray*}
 m(\bftheta,\phi) &=& E \int e^{-\bftheta. \bfx-\phi t} \eta_{\emptyset}(d\bfx)\xi_{\emptyset}(dt) \\
 &=& (1-p) \int_0^1 e^{\theta_1\log{u}+\phi\log{u}} + e^{\theta_1\log{(1-u)}+\phi\log{(1-u)}} du \\
& & \qquad +  p \int_0^1 e^{\theta_2\log{u}+\phi\log{u}} + e^{\theta_2\log{(1-u)}+\phi\log{(1-u)}} du \\
  &=& \frac{2(1-p)}{1+\theta_1+\phi} + \frac{2p}{1+\theta_2+\phi},
\end{eqnarray*}
which is finite for $(\bftheta,\phi) \in \{ \theta_1,\theta_2,\phi: \theta_1+\phi>-1, \theta_2+\phi>-1\}$.

Then we can find $\alpha(\bftheta)$ as
\begin{eqnarray*}
\alpha(\bftheta) &=& \inf\{\phi: m(\bftheta,\phi)\leq 1\} \\
&=& \inf\{\phi: \phi^2 + (\theta_1+\theta_2)\phi +\theta_1\theta_2+(2p-1)(\theta_2-\theta_1)-1 \geq 0 \}. 
\end{eqnarray*}
Thus by continuity we can solve the quadratic and take the smallest root over the set where 
$m(\bftheta,\phi)$ is defined to get
\[ \alpha(\bftheta) = -\frac12(\theta_1+\theta_2) + \frac12 \sqrt{(\theta_1-\theta_2)^2+4(2p-1)
(\theta_1-\theta_2)+4}. \]
Note that $\alpha(\bftheta)$ exists for all $\bftheta \in \br^2$ and that it is positive for all 
$\theta_1+\theta_2<0$.

We now wish to determine the Legendre transform
\[ \alpha^*(\bfa) = \inf_{\bftheta}\{ \bfa.\bftheta + \alpha(\bftheta)\}. \]
At this point it is useful to change variables and set $\varphi=\theta_1+\theta_2$ 
and $\psi=\theta_1-\theta_2$ with $\varphi,\psi\in\br$,
\[ \theta_1 = \frac{\varphi+\psi}{2}, \;\; \theta_2 = \frac{\varphi-\psi}{2}. \]
Thus our function to minimize becomes
\begin{eqnarray*}
\bfa.\bftheta + \alpha(\bftheta) 
&=& \frac12 \varphi(a_1+a_2-1) + \frac12 (a_1-a_2)\psi + \frac12 \sqrt{\psi^2+4(2p-1)\psi+4}.
\end{eqnarray*}
Hence, if $a_1+a_2-1\neq 0$, the infimum is $-\infty$ by taking $\varphi\to \pm \infty$. 
If we have $a_1+a_2-1=0$, then we are left with an equation in $\psi$ 
which we can minimise to get
\[ \psi = -2(2p-1)+\frac{4\sqrt{p(1-p)}|a_1-a_2|}{\sqrt{1-(a_1-a_2)^2}}, \]
and hence
\begin{equation}\label{1708281411}
\alpha_p^*(\bfa) = \left\{ \begin{array}{ll} 2\sqrt{p(1-p)}\sqrt{1-(a_1-a_2)^2} + (a_1-a_2)(2p-1), &
 a_1+a_2=1 \\ -\infty & a_1+a_2 \neq 1. \end{array} \right.
\end{equation}
In the case where $p=1/2$ this is
\[ \alpha^*(a_1,a_2) = \left\{ \begin{array}{ll} \sqrt{1-(a_1-a_2)^2}, & a_1+a_2 = 1 \\
-\infty & a_1+a_2\neq 1. \end{array} \right. 
\]
That is, for $0\leq a_1 \leq 1$ we have
\[ \alpha^*(a_1,1-a_1) = 2\sqrt{a_1(1-a_1)} = 2\sqrt{a_1a_2} \]
and $\alpha^*(a_1,a_2)=-\infty$ for all other values of $\bfa=(a_1,a_2)$.

Unfortunately, as $\alpha^*(\bfa)$ is only finite on a closed set, we cannot apply the original theorem 
of Biggins but give our own version for this degenerate set up. 

\begin{thm}\label{1708261300}
Suppose that $A$ is a closed convex set in $\br_+^2$ with a non-empty interior.
Let $\beta = \sup\{\alpha^*(a): a\in A\}$.
\begin{enumerate}
\item If $A \cap\{(x,y):x+y=1\}=\emptyset$, then
\[ t^{-1} \log{N_t(tA)} \to -\infty, \;\; t\to\infty, \;\; a.s. \]
\item If $A \cap\{(x,y):x+y=1\}\neq \emptyset$, then setting $\beta = \sup\{\alpha^*(a):a\in A\}$ we have 
\[ t^{-1}\log N_t(tA) \to \beta, \;\; t\to\infty\;\; a.s. \]
\end{enumerate}
\end{thm}

\emph{Proof:}
In order to prove this theorem we use the one-dimensional version of the result from \cite{Big95} 
and follow the same line of reasoning as \cite{Big94}. At this point we will just prove the version 
where $p=1/2$. The extension to general $p$ is a minor modification.

We consider the projections of the model onto lines through the origin at angle $\psi$. That is the points
$(0,y) \to y\sin \psi$ and $(x,0) \to x\cos\psi$. We now consider the one-dimensional projected general 
branching random walk which has offspring distribution as
\[ \begin{array}{ll} (-\cos\psi\log{U}, -\cos\psi\log{(1-U)}) & \mbox{ with probability } 1/2 \\
 (-\sin\psi\log{U}, -\sin\psi\log{(1-U)}) & \mbox{ with probability } 1/2. \end{array} \]
 
The Theorem of Biggins \cite{Big95} will then apply to this case for all points except possibly 
at $\tilde{a}_{\psi}:=\sup\{a: \alpha^*_{\psi} (a)>-\infty\}$, the right boundary of finiteness of $\alpha^*$.

Let $n_{\psi} = (\cos(\psi),\sin(\psi))$ be a unit vector with angle $\psi$.
By construction the moment generating function for the projection of the GBRW onto the line in the 
direction $n_{\psi}$ is given by
\[ m_{\psi}(\theta,\phi) = m(\theta n_{\psi},\phi) = \frac{1}{1+\cos(\psi) \theta+\phi} +
\frac{1}{1+\sin(\psi)\theta+\phi} \]
and similarly we have
\[ \alpha_{\psi}(\theta) = \alpha(\theta n_{\psi}) = -\theta c_{\psi} + \sqrt{\theta^2 d_{\psi}^2 +1},\]
where
\[ c_{\psi} = \frac{\cos(\psi)+\sin(\psi)}{2}, \;\; d_{\psi} = \frac{\cos(\psi)-\sin(\psi)}{2}. \]
Computing the Legendre transform of $\alpha_{\psi}(\theta)$ as before we obtain, assuming $\psi \neq \pi/4$,
\[
\alpha_{\psi}^*(a) = \left\{ \begin{array}{ll} \sqrt{1-\left(\frac{c_{\psi}-a}{d_{\psi}}\right)^2}, 
& |c_{\psi}-a|\leq d_{\psi}, \\
-\infty, & |c_{\psi}-a|>d_{\psi}. \end{array} \right. \]
Rewriting this we have
\[ \alpha_{\psi}^*(a) = \left\{ \begin{array}{ll} 2\sqrt{\frac{(\cos(\psi)-a)(a-\sin(\psi))}
{(\cos(\psi)-\sin(\psi))^2}}, & \sin(\psi)\wedge \cos(\psi) \leq a \leq \cos(\psi)\vee \sin(\psi), \\
-\infty, & a < \sin(\psi)\wedge\cos(\psi), \quad a > \cos(\psi)\vee\sin(\psi). \end{array} \right. \]
In the case where $\psi=\pi/4$ we see that $\alpha_{\pi/4}(\theta) = -\frac{\theta}{\sqrt{2}} +1$ and thus 
\[ \alpha^*_{\pi/4}(a) = \left\{ \begin{array}{ll}1, &  a=1/\sqrt{2}, \\ 0, & a \neq 1/\sqrt{2}. 
\end{array} \right. \] 
Let 
\[ \cL_{c} =\bigcap_{\psi\in [0,2\pi)} \{\bfa: \alpha^*_{\psi}(\bfa.n_{\psi})\geq c\}, \]
where
\[ \alpha^*_{\psi}(\bfa.n_{\psi}) =  2\sqrt{\frac{((1-a_1)\cos(\psi)-a_2\sin(\psi))(a_1\cos(\psi)-
(1-a_2)\sin(\psi))}{(\cos(\psi)-\sin(\psi))^2}}.\]
It is easy to see that the line segment $\{\bfa: 2\sqrt{a_1 a_2}\geq c, a_1+a_2=1\}\subset 
\{\bfa: \alpha^*_{\psi}(\bfa.n_{\psi})\}$ for all $\psi\in [0,\pi/4)\cup(\pi/4,\pi/2]$.
We also observe that as $\psi\to\pi/4$ we have for $|a_1+a_2-1+\epsilon(a_1-a_2) +o(\epsilon)|<\epsilon$
\[ \alpha^*_{\pi/4-\epsilon}(\bfa.n_{\pi/4-\epsilon}) = \frac{1}{\epsilon}\sqrt{(a_1+a_2-1 + \epsilon(a_1-a_2+1))
(1-a_1-a_2+\epsilon(1-a_1+a_2))} + O(\epsilon).\]
Hence we see that at $\pi/4$, by letting $\epsilon\to 0$, we must have
\[ \alpha^*_{\pi/4}(\bfa.n_{\pi/4})=\sqrt{1-(a_1-a_2)^2} \mbox{ if }  a_1+a_2=1, |a_1-a_2|<1. \]
That is 
\[ \alpha^*_{\pi/4}(\bfa.n_{\pi/4})=2\sqrt{a_1(1-a_2)} \mbox{ if }  a_1+a_2=1, 0< a_1 <1. \]
If $a_1+a_2 \neq 1$, then $ \alpha^*_{\pi/4}(\bfa.n_{\pi/4})=-\infty$. 

We now remark on the boundary. We know that the one-dimensional results of Biggins hold for all $a$ 
except possibly at $\tilde{a}=\sup\{a: \alpha_{\psi}^*(a)>-\infty\}$. In the case of the point $\tilde{a}$ 
we know that $\exp(\alpha^*_{\psi}(\tilde{a})t)$ is an upper bound for
the number of particles at $\tilde{a} t$ and beyond. In our setting we have, writing 
$\tilde{a}_{\psi}=\sup\{a: \alpha_{\psi}^*(a)>-\infty\}$, that $\alpha^*_{\psi}(\tilde{a}_{\psi}) = 0$ 
for all $\psi\neq \pi/4 \in [0,\pi/2]$. At the point $\psi=\pi/4$ we have $\tilde{a}_{\pi/4}=1/\sqrt{2}$ 
and $\alpha^*_{\pi/4}(1/\sqrt{2})=1$. It is easy to check in this case that by construction we have with 
probability one a particle born at position $-(\log U)/\sqrt{2}$ at time $-\log U$ and one at 
$(-\log(1-U))/\sqrt{2}$ at time $-\log(1-U)$ and hence at time $t$ all $e^t$ particles are close to 
$t/\sqrt{2}$ and we have the exponential growth at rate 1. For the other cases we note that as 
$\alpha^*_{\psi}(\tilde{a})=0$ 
we are considering the behaviour of the rightmost particle $R^{\psi}_t$. By using \cite{Big95}~Corollary~2 
we know that, under the conditions of our theorem, that $R^{\psi}_t/t \to \gamma$, where 
$\gamma:=\inf\{a: \alpha_{\psi}^*(a)<0\}$. Thus we have $\gamma=\tilde{a}$ and at this value
\[ \frac{\log N_t(t\gamma,\infty)}{t} \to 0. \]
Thus the limit result for the number of particles in $(t a,\infty)$ holds at $a=\tilde{a}$ too.

Hence $\cL_c$ will be equivalent to the line segment in $\mathbb{R}^2_+$ where $\alpha^*(\bfa)$ is above 
$c$. 

The claims of the theorem now follow from a minor modification of the proof of Theorem~4.1 of 
\cite{Big94}.

\enpf

We now treat the original problem with the Theorem we have just obtained. For a set $S\subset [0,1]^2$ we 
ask about $\tilde{N}_t(S)$, the number of rectangles with side lengths $(a,b)$ in $S$ at time $t$. We 
transform the function $\alpha^*$ to $\gamma^*$ as
\[ \gamma^*(a,b) = 2\sqrt{\log{a}\log{b}}, \;\; ab=e^{-1}, \]
and $\gamma^*=0$ for all other values of $a,b$. 

\begin{cor}\label{cor:gbrw}
Let $S$ be a closed convex set with non-empty interior in $[0,1]^2$. 
\begin{enumerate}
\item If $S\cap\{(a,b): ab=e^{-1}\}=\emptyset$, then
\[ t^{-1} \log{\tilde{N}_t(e^{-t}S)} \to -\infty, \;\; t\to\infty, \;\; a.s. \]
\item If $S \cap\{(a,b): ab=e^{-1}\}\neq \emptyset$, then setting $\beta = \sup\{\gamma^*(a,b):(a,b)\in S\}$ 
we have 
\[ t^{-1}\log{\tilde{N}}_t(e^{-t}S) \to \beta, \;\; t\to\infty\;\; a.s. \]
\end{enumerate}
\end{cor}

This can be viewed as roughly giving for instance that if $S_{a,b}$ is a closed convex set with non-empty 
interior such that $S_{a,b} \cap \{(x,y):xy=e^{-t}\} = \{(a,b)\}$, then for large $t$
\[ \tilde{N}_t(S_{a,b}) \approx \exp(2\sqrt{-\log{a}(t+\log{a}}). \]


The quantity we are interested in is the density of interfaces of a given length. At this stage we have 
a description of all the rectangles in the unit square after a time $t$. The interfaces come from the 
lengths of the sides of these rectangles. That is, when a rectangle is split, aside from two new rectangles 
we have an interface whose length is the length of the side which is not split. Thus, if we consider
$F_t(x)$ to be the number of interfaces that are greater than $x$ at time $t$, then we can write
this in terms of the branching process as 
\[ F_t(x) = \sum_{\bfi \in {\cal T}} I_{\{(Z_t^{\bfi})_i\leq -\log{x}, i=1 \mbox{ or }2\}}
I_{\{\sigma_{\bfi}<t\}}. \]
Hence we are counting the branching process with a 0-1 characteristic and can apply our Theorem to 
conclude that
\[ \lim_{t\to\infty}\frac{\log{F_t(tx)}}{t} = 2\sqrt{-\log{x}(1+\log{x})}. \]

\subsection{Splitting independently}\label{sec:indep}

Our calculations so far supposed that we always split the largest rectangles first. 
An alternative is to choose to split each rectangle independently of the rest. Our approach 
via a general branching random walk process allows us to use any distribution for the splitting time 
but we can obtain more precise limit theorems in the case where the split is at constant rate, in 
which case we have a Markov branching random walk. In this setting small pieces divide in the same 
way as the large pieces and this will give, in general, a greater disparity in the sizes of the 
individual rectangles at a given time. 

We will compute $\alpha^*(\bfa)$, the function which determines
the asymptotic growth rate in the general branching random walk for this version of the model.
In general, if the birth time process $\xi$ is independent of the spatial displacement $\eta$, we will have
\[ m(\bftheta,\phi) = E \int_{\mathbb{R}^d} e^{-\bftheta. \bfx} \eta_{\emptyset}(d\bfx) E\int_0^{\infty} 
e^{-\phi t} \xi_{\emptyset}(dt), \;\; \bftheta\in \R^d, 
\phi \in \R_+. \]
Let $\hat{\xi}(\phi) = E \int_0^{\infty} \exp(-\phi t)\xi(dt)$ so that
\[ m(\bftheta,\phi) = \hat{\xi}(\phi)\left(\frac{1}{1+\theta_1}+\frac{1}{1+\theta_2}\right) \]
where we assume that $m(\bftheta,\phi)<\infty$ in a neighbourhood of the origin.
Now
\[ \alpha(\bftheta) = \hat{\xi}^{-1}\left(\left(\frac{1}{1+\theta_1}+\frac{1}{1+\theta_2}\right)^{-1}\right), 
\]
and $\alpha^*$ will be the Legendre transform defined by
\[ \alpha^*(\bfa) = \inf_{\bftheta}\{ \bfa.\bftheta + \alpha(\bftheta)\}. \]
To compute this we write $\psi_i = 1+\theta_i$ and $b(\psi_1,\psi_2)= 1/\psi_1 + 1/\psi_2$. Now, 
expressing the problem in terms of the variables $\psi_i$, we have the first order conditions
\begin{eqnarray*}
a_1 -\frac{1}{b^2} \frac{\partial b}{\partial \psi_1} (\hat{\xi}^{-1})'(\frac{1}{b}) &=& 0 \\
a_2 -\frac{1}{b^2} \frac{\partial b}{\partial \psi_2} (\hat{\xi}^{-1})'(\frac{1}{b}) &=& 0. 
\end{eqnarray*}
Thus we must have
\[ a_2 \frac{\partial b}{\partial \psi_1} = a_1 \frac{\partial b}{\partial \psi_2}, \]
that is $\psi_1= \sqrt{a_2/a_1} \psi_2$. Using this to express our first order equations in terms of 
$\psi_2$ we have that $\psi_2$ must satisfy
\[ a^2 = -(\hat{\xi}^{-1})'(\frac{\psi_2\sqrt{a_2}}{a}), \]
where $a = \sqrt{a_1}+\sqrt{a_2}$. Using this in our expression for $\alpha^*$ we have, writing 
$\zeta(x) = [\hat{\xi}'(x)]^{-1}$ for the inverse function of the negative of the derivative of 
the Laplace transform in time, that $\psi_2 = a\hat{\xi}(\zeta(-1/a^2))/\sqrt{a_2}$ and hence
\[ \alpha^*(\bfa) = a^2\hat{\xi}(\zeta(-\frac{1}{a^2}))-a_1-a_2+\zeta(-\frac{1}{a^2}). \]

In the case where there is bias in the offspring distribution, with $p$, the probability of a 
horizontal edge, then similar calculations give the same formula where now $a = \sqrt{2pa_1} + 
\sqrt{2(1-p)a_2}$.

We note that the same argument allows us to handle the three dimensional case and we find that for 
$\bfa = (a_1,a_2,a_3)$ 
\[ \alpha^*(\bfa) = a^2\hat{\xi}(\zeta(-\frac{1}{a^2}))-a_1-a_2-a_3+\zeta(-\frac{1}{a^2}), \]
where now $a = \sqrt{a_1}+\sqrt{a_2}+\sqrt{a_3}$ and $\hat{\xi}$ and $\zeta$ are as before. Again 
adding bias just changes the expression for $a$ to $a= \sqrt{3p_1a_1}+\sqrt{3p_2a_2}+\sqrt{3p_3a_3}$, 
where $p_1,p_2,p_3$ with $p_1+p_2+p_3=1$ are the probabilities for the three axial directions.

We now consider a couple of natural cases.

\subsubsection{Constant rate case}

The simplest version is the case where individuals reproduce independently at a constant rate. The GBRW 
becomes a Markov branching random walk and we can apply more precise limit theorems. As it is still a 
GBRW we stay within that framework initially to compute the asymptotic growth rate as before. 
Firstly the Laplace transform of the offspring birth and location measure is 
\[ m(\bftheta,\phi) = \left(\frac1{1+\theta_1} +\frac{1}{1+\theta_2}\right)\frac{1}{1+\phi}. \]
Hence
\[ \alpha(\bftheta) = \frac1{1+\theta_1} +\frac{1}{1+\theta_2} -1. \]
Inverting this by setting
\begin{equation}
\theta_1= a_1^{-1/2}-1, \;\; \theta_2 = a_2^{-1/2}-1 \label{eq:theta12cr}
\end{equation}
gives
\begin{eqnarray}\label{1708271448}
 \alpha^*(\bfa) = 1-(1-\sqrt{a_1})^2-(1-\sqrt{a_2})^2, \;\; a_1,a_2>0. 
\end{eqnarray}
As $\alpha^*(\bfa)>0$ in an open set we can apply a transformed version of Theorem~\ref{thm:gbrw}. We let 
\begin{equation}
 \gamma^*(a,b) = \alpha^*(-\log{a},-\log{b}) = 1-(1-\sqrt{-\log{a}})^2-(1-\sqrt{-\log{b}})^2, \;\; 0< a,b<1 
 \label{eq:gamma-indep}
\end{equation}
and $S$ be a closed convex set in $[0,1]^2$. Then for $\beta = \sup_{(a,b)\in S} \gamma^*(a,b)$ we have 
if $\beta>0$ 
\[ \lim_{t\to \infty} \frac{1}{t}\log{\tilde{N}_t(e^{-t}S)} = \beta, \;\; a.s. \]
while if $\beta<0$, then for any $\delta>\beta$
\[ \lim_{t\to\infty} e^{-\delta t}\tilde{N}_t(e^{-t}S) = 0, \;\; a.s. \]

In \cite{Uch} it is shown that for a Markov branching random walk, where the individual gives birth 
at the moment of their death, that there is a finer limit theorem on the growth of the number of 
individuals in a certain set. The result states, in our notation, that
for a fixed $\bfa$ such that $\alpha^*(\bfa)>0$, we have
\[ \lim_{t\to\infty} t\exp(-\alpha^*(\bfa)t) N_t(t\bfa + D) = \frac{a_1^{3/4} a_2^{3/4}}{\pi} \int_D 
e^{\theta_1x+\theta_2y} dxdy W_{\bfa}, \;\; a.s. \]
where $W_{\bfa}>0$ is a non-degenerate mean one random variable and $\theta_1,\theta_2$ are as in \eqref{eq:theta12cr}.

\subsubsection{Discrete generations}\label{201809141440}

An alternative
version of this model is the one in which we work in generations, in that each rectangle splits at a 
fixed time 1. This can be encoded by a classical branching random walk and a realization of the associated
microstructure is shown in Figure~\ref{2006172133}. 
Fitting this into the framework of this section we have $\xi = \delta_1$ and $\hat{\xi}(\phi)=
\exp(-\phi)$. Applying the same calculations we have
\[ \alpha(\bftheta) = \log\left(\frac{1}{1+\theta_1}+\frac{1}{1+\theta_2}\right), \]
and hence
\begin{eqnarray}\label{1708262241}
 \alpha^*(\bfa) = 1-a_1-a_2+2\log(\sqrt{a_1}+\sqrt{a_2}). 
\end{eqnarray}
The constant rate model can be approximated by a geometric random variable with a probability $\Delta t$ 
of success. We consider the case where 
$\Delta t \leq 1$ and define the life time distribution of the rectangle as
\[ \xi = \sum_{n=1}^{\infty} \Delta t(1-\Delta t)^{n-1} \delta_{n\Delta t}. \]
This corresponds to a model in which each rectangle splits with probability $\Delta t$ at each 
discrete time step $n\Delta t, \;n\in \mathbb{N}$. Computing the Laplace transform for the life 
time distribution we have 
\[ \hat{\xi}(\phi) = \frac{\Delta t \exp(-\phi \Delta t)}{1-(1-\Delta t)\exp(-\phi \Delta t)}, \]
and hence, keeping track of the dependence on $\Delta t$, we have
\[ \alpha_{\Delta t}(\bftheta) = \frac{1}{\Delta t} \log\left(1-\Delta t + \Delta t\left(\frac{1}
{1+\theta_1}+\frac{1}{1+\theta_2}\right)\right). \]
Computing the Legendre transform, and writing $a=\sqrt{a_1}+\sqrt{a_2}$, we have
\begin{eqnarray*}
\alpha_{\Delta t}^*(\bfa) &=& a \frac{\sqrt{a^2(\Delta t)^2+ 4(1-\Delta t)}- 
a\Delta t}{2(1-\Delta t)}  \\ & & \qquad + \frac{1}{\Delta t}\log\left(1+\Delta t\left(\frac{2a(1-\Delta t)}
{\sqrt{a^2(\Delta t)^2 + 4(1-\Delta t)}-a\Delta t} -1 \right)\right)-a_1-a_2. 
\end{eqnarray*}

It is possible to check that in the limit we recover the constant rate case by letting $\Delta t \to 0$ 
and we get the fixed time version by letting $\Delta t \to 1$. In the numerical work we take $\Delta t$ 
small and it is not hard to see from the asymptotics that for $\Delta t$ small
\begin{eqnarray*}
 \alpha_{\Delta t}^*(\bfa) &= & 2(\sqrt{a_1}+\sqrt{a_2}) -a_1-a_2-1 + o(\Delta t), \\
 &=& \alpha^*_c(\bfa) + o(\Delta t),
\end{eqnarray*}
where $\alpha^*_c(\bfa)$ is the Legendre transform given in \eqref{1708271448}.

%
%

%

\section{An alternative view and a power law}\label{1708081205}

We will now consider the lengths of horizontal interfaces in a complete fragmentation of the unit square.
Let $N^p_h(x)$ denote the number of horizontal interfaces in the unit square 
which are of length greater than or equal to $x$ when the probability of horizontal edges is $p$.
It is a simple observation that if $p>1/2$, this quantity will be infinite with 
positive probability, as there is a positive probability that the number of edges of unit length goes 
to infinity, as the number of edges of unit length is a supercritical Galton-Watson branching process. 
In the case where $p=1/2$ the random variable $N^p_h(x)$ has infinite mean and we will not discuss it further.
As a result we will only consider the case $p<1/2$ and will see that $N^p_h(x)$ is almost surely finite 
and will give a limit theorem for its scaling with $x$.

In order to do this we introduce another general branching process. 
As we are tracking the number of horizontal edges of a certain size, if a rectangle is split in the horizontal 
direction, that will correspond to having a horizontal edge of exactly the same length as the parent and two
new rectangles which could produce more horizontal edges of the same size. A split in the vertical direction 
will lead to new rectangles which will produce edges of smaller length. To map this into a general branching 
process we will regard the horizontal axis as time and an individual's offspring will correspond to two 
sets of edges; those that are the offspring of the first child which all have the same horizontal length and 
those that are the offspring of the second child which all have the same horizontal length. As $p<1/2$, the 
number of horizontal edges generated at a single location is a subcritical branching 
process and we will see that it is straight forward to compute its law.

We now consider the family of individuals corresponding to horizontal edges at a given location. There are 
two types of edge; those that have given rise
to further horizontal edges of the same size and those that just yield edges of a smaller size. 
It is this second type that constitute
the children for our general branching process, as they will continue to reproduce. The first type just 
need to be counted as they contribute to the number of interfaces. 
Thus we will have a general branching process where the offspring law is $\tilde{X}_1$ individuals at time 
$-\log U$ and $\tilde{X}_2$ at time $-\log(1-U)$, where $\tilde{X}_1,\tilde{X}_2$ are independent random variables
determined by the law of the edges that reproduce. 
As our initial condition in the model is a unit square the initial condition
for the general branching process is not necessarily a single individual. It will be the reproducing members 
of the subcritical branching process located at the origin. 

We now compute the law of the sub-crititcal branching process.
Let $(S_n)_{n=0}^{\infty}$ be a simple random walk on the positive integers with 
\[ P(S_n-S_{n-1}=1) = p, \;\; P(S_n-S_{n-1}=-1) = 1-p, \]
and let $T_0=\inf\{n: S_n=0\}$.

\begin{lem}
The law of the number of offspring at a given $x$ is that of $T_0$ given that $S_0=1$. 
\end{lem}

{\it Proof:}
We observe that if a parent with horizontal length $x$ has two children with the same horizontal length 
(with probability $p$) the number of horizontal components of length $x$ increases by one. If the two
children are of horizontal length less than $x$ (with probability  $1-p$), then the number of horizontal 
components of size $x$ decreases by one. In this way we see that the number of horizontal components of 
a given size $x$ behaves as a simple random walk with drift. At the time $T_0$ there are no more horizontal 
edges of the parent length that can be produced.
\enpf

With this observation we can construct the law of our general branching process. As observed the number of 
horizontal components behaves like a simple random walk. We observe that each upstep of the walk, which
corresponds to the increase in the number of components does not lead to a `birth' in the future. The down 
steps correspond to components that will reproduce in the future. Thus we can observe that there are
$(T_0+1)/2$ offspring that will reproduce in the same way as the parent. There are also $(T_0-1)/2$ 
individuals that correspond to horizontal edges of length $-\log{x}$. 
It is these pieces that we wish to count and so we use the following general branching process.

Let $(\xi, L , \chi)$ be the triple consisting of the birth point process $\xi$, the lifelength $L$ and 
the characteristic counting function $\chi$. All are defined using a $U[0,1]$-random variable $U$ and 
$T_0,\tilde{T}_0$, two independent copies of $T_0$. The birth point process is 
\[ \xi(dt) = \frac{T_0+1}{2}\delta_{-\log U}(dt) + \frac{\tilde{T}_0+1}{2}\delta_{-\log{(1-U)}}(dt). \]
The lifelength is $\max\{-\log{U},-\log{(1-U)}\}$, that is the parent dies at the birth of its final child.
The counting function $\chi$ is used to capture the number of edges of length at least 
$e^{-t}$ produced by an individual and is thus given by
\[ \chi(t) = \frac{T_0-1}{2}I_{\{t\geq -\log{U}\}} + \frac{\tilde{T}_0-1}{2}I_{\{t\geq -\log{(1-U)}\}}. \]
We write 
\[ Z^{\chi}(t) = \sum_{\bfi \in \Sigma} \chi_{\bfi}(t-\sigma_{\bfi}), \]
so that $Z^{\chi}(t)$ counts all the edges of length greater than $e^{-t}$. Theorem~\ref{thm:sllngbp} shows that, 
provided there are not too many offspring too far in the future, there will be a strong law of large numbers
for this process. We first find the Malthusian parameter $\alpha$ as the solution to
\begin{eqnarray*}
1 &=& \be \sum_{i=1}^{\xi(\infty)} e^{-\alpha\sigma_i} \\
&=& \be\left( U^{\alpha} \frac{T_0+1}{2}\right) + \be \left((1-U)^{\alpha} \frac{\tilde{T}_0+1}{2} \right) \\
&=& \frac{\be T_0+1}{\alpha+1},
\end{eqnarray*} 
using the independence of $U$ and $T_0, \tilde{T}_0$.
That is the Malthusian parameter is given by $\alpha = \be T_0 = 1/(1-2p)$.

Recalling the notation from the section on the general branching process we then set
\[ z^{\chi}(t) = e^{-\alpha t} \be Z^{\chi}(t), \]
and
\begin{eqnarray*}
 u^{\chi}(t) &=& e^{-\alpha t} \be \chi(t) \\
 &=& e^{-\alpha t} \be \left( \frac{T_0-1}{2}I_{\{t\geq -\log{U}\}} + \frac{\tilde{T}_0-1}{2}
 I_{\{t\geq -\log{(1-U)}\}} \right) \\
 &=& e^{-\alpha t} \be\left( (T_0-1) I_{\{t\geq -\log{U}\}} \right) \\
 &=& (\alpha-1) e^{-\alpha t} \bp(U\geq e^{-t}) \\
 &=& (\alpha-1) e^{-\alpha t}(1-e^{-t}).
 \end{eqnarray*}
There is also a fundamental martingale $M$ determined by the weighted birth times of the coming generation:
\[ M_t = \sum_{j:\sigma_i\leq t< \sigma_{ij}} e^{-\sigma_{ij}}. \]
 
\begin{thm}\label{thm:lim1}
For the counting process $Z^{\chi}(t)$ we have
\[ \lim_{t\to \infty} e^{-\alpha t} Z^{\chi}(t) = (1-\frac{1}{\alpha})M_{\infty} = 2p M_{\infty}, \;\; a.s. \]
where $M_{\infty} = \lim_{t\to\infty} M_t$ exists and is non-degenerate.
\end{thm} 
 
{\it Proof:} 
In order to apply Theorem~\ref{thm:sllngbp}, the strong law of large numbers for the general branching 
process, we need to check the conditions. Firstly observe that the expected family size is finite and 
the measure $\nu_{\alpha}$ is non-lattice. The characteristic $\chi$ is increasing in $t$ with
$\be\chi(\infty) = \alpha-1<\infty$ and hence, taking the function $h(t)=\exp(-\epsilon t)$ for $\epsilon<\alpha$, in the conditions of Theorem~\ref{thm:sllngbp}, we have
\[ \be \left( \sup_{t \geq 0} e^{-(\alpha-\epsilon) t} \chi(t) \right) \leq \be\chi(\infty) = \alpha-1<\infty. \]
Thus we can apply the Theorem and all we need is to compute the limit in the renewal equation. That is
\[ z^{\chi}(t) \to \frac1{\mu_1}\int_0^{\infty} u^{\chi}(t) dt, \]
where $\mu_1 = \be \int_0^{\infty} t e^{-\alpha t} \xi(dt)$. Straightforward calculations give
\begin{eqnarray*}
\int_0^{\infty} u^{\chi}(t) dt &=& \int_0^{\infty}  (\alpha-1) e^{-\alpha t}(1-e^{-t}) dt\\
&=& \frac{\alpha-1}{\alpha(\alpha+1)},
\end{eqnarray*}
and
\begin{eqnarray*}
\mu_1 &=&  \be \int_0^{\infty} t e^{-\alpha t} \xi(t) dt \\
&=& \be \left( (-\log{U}) U^{\alpha} \frac{T_0+1}{2} +(- \log{(1-U)})(1-U)^{\alpha}\frac{\tilde{T}_0+1}{2} 
\right) \\
&=& (\alpha+1)\be\left( (-\log{U}) U^{\alpha} \right)\\
&=& \frac{1}{\alpha+1}.
\end{eqnarray*}
Putting these two calculations together we have
\[ z^{\chi}(t) \to 1-1/\alpha = 2p, \]
as required.
\enpf

\begin{rem}{\rm 
An alternative way to see the growth rate is to consider the general branching process in which individuals 
can give birth immediately, so that for individual $\bfi$ we have $\xi_{\bfi}(dt) = 2\delta_{t=0}I_{B_{\bfi}=1} 
+ (\delta_{t=-\log{U}} + \delta_{t=-\log{(1-U)}})I_{B_{\bfi}=0}$. It is then straightforward to
calculate $\alpha$ such that $1=\be \sum_{i=1}^{\xi(\infty)} \exp(-\alpha \sigma_i)$ and arrive at
the same result. We do not follow this route as such offspring distributions are not usually part of 
general branching process theory.}
\end{rem}

As a corollary we have our power law limit theorem.
%

\begin{thm}\label{1708280114}
There exists a strictly positive random variable $W = \sum_{i=1}^{(T_0+1)/2} M_{\infty}^{(i)}$, where 
$M^{(i)}_{\infty}$ are independent copies of $M_{\infty}$,  such that
\[ \lim_{x\to 0} x^{1/(1-2p)}N^p_h(x) = 2p W, \;\; a.s. \]
\end{thm}

{\it Proof:} We observe that the initial condition for the process is the number of reproducing individuals 
at time 0 which is given by $(T_0+1)/2$. We also have $(T_0-1)/2$ non-reproducing individuals which are 
counted. We will write $Z^{\chi}_i(t)$ for the evolution of the general branching process started from the 
initial individual $i$. Thus, making the time transformation that $t=-\log{x}$, we have
\[ N^p_h(x) = \frac{T_0-1}{2} + \sum_{i=1}^{(T_0+1)/2} Z^{\chi}_i(-\log{x}). \]
Hence
\[ \lim_{x\to 0} x^{1/(1-2p)} N^p_h(x) = \lim_{t\to\infty} \left( e^{-t/(1-2p)}\frac{T_0-1}{2} + 
\sum_{i=1}^{(T_0+1)/2} e^{-t/(1-2p)}Z^{\chi}_i(t) \right). \]
The result then follows from Theorem~\ref{thm:lim1}.
\enpf

\begin{rem}{\rm 
1). The three dimensional case of the model enables us to compute the exponent without requiring bias and 
we discuss this in the next section along with a range of other variations on the basic model.

2). The power law result is a function of the whole fragmentation and the time parameter used to generate
the fragmentation has no impact on this.}
\end{rem}

\section{Alternative models}\label{1801071844}

We give brief accounts of how our methods can be extended to alternative models. In particular we use a 
GBRW analysis and the power law analysis via the GBP and label the approaches by (G) and (P) respectively.
In the case (G) we will compute the function $\alpha^*(\bfa)$ for the branching process. This is related
to the decay function for rectangles via the transformation $\gamma^*(a,b) = \alpha^*(-\log{a},-\log{b})$.

\subsection{The three dimensional case}\label{1708031736}

Here we take the unit cube and split it along the three axes using the same mechanism as in two dimensions. 

(G) The three dimensional problem leads to a GBRW in the same way, although now we have three directions 
in which the individuals can move. The moment generating function is given by
\[ m(\bftheta,\phi) = \frac{2}{3}\left(\frac{1}{1+\theta_1+\phi}+\frac{1}{1+\theta_2+\phi} +
 \frac{1}{1+\theta_3+\phi}\right) \]
and this can be used to find $\alpha(\bftheta)$ as the solution to a cubic equation.
The computation of $\alpha^*$ is then best done numerically.

%

(P) In the alternative power law view we obtain Theorem~\ref{thm:3dP}. If we take our two dimensional 
picture of the previous case but now take as the time axis the area of the horizontal slabs we see 
that a horizontal slab preserves its area, if there is a horizontal slice, and we produce an 
interface with that area. When we split in a different direction the horizontal area is uniformly 
split into two pieces. Thus we have exactly the same GBP as before but now the critical probability 
is 1/3 as this is the chance of a horizontal split. Directly applying our previous result we see 
that if $N_h(a)$ is the number of horizontal slabs of area $a$, then 
\[ \lim_{a\to 0} N_h(a) a^{3} = \frac23 W,\;\; P-a.s. \]
Thus we have Theorem~\ref{thm:3dP}.

\subsection{Splitting right triangles}

In the general case of the martensitic phase transition there are a set of compatible angles which determine
the directions of the plates that form. These cannot in general be mapped to axial directions and the 
different component types cannot be kept track of using such a simple coordinate system as we use here. We discuss a simple modification of our model to illustrate the difficulties with even the simplest modifications.

We start with a right angled triangle and then split it into four right angled triangles by choosing a 
point in the interior according to the uniform distribution and then drawing a line parallel to one of the right angled edges until it hits the hypotenuse and then joining this point to the other edge. 
This leaves two triangles and a rectangle which can be further split into two triangles by a line along the 
main diagonal. The distribution of the side lengths of each triangle can be computed.

(G) In this case we have,
\[ (a,b) \to \left\{\begin{array}{ll} \{(aU,bU), (a(1-U),bU), (a(1-U),bU),(a(1-U),b(1-U)) \} & 
\mbox{ probability } 1/2 \\  
\{(aU,bU), (aU,b(1-U)),(aU,b(1-U)),(a(1-U),b(1-U)) \} & \mbox{ probability } 1/2 \end{array} \right. \]
with birth times given by their size, that is $U^2, U(1-U), (1-U)^2$. Again we take logs to map this into a
general branching random walk. The Laplace transform for this GBRW is given by
\[ m(\bftheta,\phi) = \frac{2}{1+\theta_1+\theta_2 + 2\phi} + 2\frac{\Gamma(\theta_1+\phi+1)
\Gamma(\theta_2+\phi+1)} {\Gamma(\theta_1+\theta_2+2\phi+2)} \]
and it does not appear possible to solve this analytically to find $\alpha(\bftheta)$.

(P) If we consider the alternative approach, in which the $x$-axis becomes time, we see that there is no chance 
of explosion as all triangles decrease in area and the interfaces produced also decrease in length. An individual 
of side length 1 has three interfaces of length $U, 1-U$ and $\sqrt{U^2+(1-U)^2}$ along with the four children. We 
can count these interfaces using a characteristic for the general branching process run with time given by the
negative logarithm of the horizontal side length. In this framework the number of particles with side length 
around $x$ is the number of individuals alive at around time $-\log{x}$. If we calculate the Malthusian 
parameter, $\gamma$, for this general branching process 
\[ 1 = 3 E U^{\gamma} + E (1-U)^{\gamma} = \frac{4}{\gamma+1}, \]
we see $\gamma = 3$. If $N(x)$ is the number of interfaces of length greater than $x$ we can 
express it in terms of $Z^{\chi}_t$, the total population in the general branching 
process up to time $t$ counted according to a characteristic $\chi$, which counts the interfaces 
of length $e^{-t}$. As $\chi$ is bounded by 3, we can apply the Theorem~\ref{thm:sllngbp} to see that
there is a constant $m^{\chi}$ such that
\[ \lim_{t\to\infty} e^{-3t} Z^{\chi}_t = m^{\chi} W, \;\; P-a.s.\]
In terms of interface length
\[ \lim_{x\to 0} N(x)x^3 = m^{\chi} W, \;\; P-a.s. \]


\begin{rem}{\rm
The move to more other angles leads to more complications. For example
imagine that we start with a regular hexagon and split it along any of the three directions 
which are parallel to its sides with equal probability in the same way as for the square. It is clear 
that the offspring of the original will consist of a hexagon, no longer regular, and a trapezium.
As we iterate then we will have irregular shapes appearing but they will all have either 3, 4, 5 
or 6 sides and can be summarized by a set of up to 6 numbers. An exact analysis looks to be very 
challenging.}
\end{rem}

\subsection{Non-uniform splitting}
\label{1809241600}
In practice nucleation sites are not equally likely to occur at any point uniformly across 
the square. In order to model this we could assume that the random variable $U$ in the model is no 
longer uniformly distributed but comes from some other distribution on $[0,1]$. We consider the 
case of the Beta distribution. 

Let $U$ have a Beta$(\alpha,\alpha)$ distribution in which the density is given by
\[ f_{\alpha}(x) = \frac{\Gamma(\alpha)^2}{\Gamma(2\alpha)} x^{\alpha-1}(1-x)^{\alpha-1}, \;\;0<x<1. \]
At $\alpha=1$ this is uniform but as we increase $\alpha$ the density function becomes more peaked around 1/2.
We note that if $U$ is Beta$(\alpha,\alpha)$ distributed then
\begin{equation}
 EU^{\gamma} =  \frac{\Gamma(\gamma+\alpha)\Gamma(2\alpha)}{\Gamma(\gamma+2\alpha)\Gamma(\alpha)}. \label{eq:bet}
\end{equation}

(G) Here we assume $p=1/2$. Using \eqref{eq:bet} in the case of Beta$(\alpha,\alpha)$, we have
\begin{equation}
 m(\bftheta,\phi) = \frac{\Gamma(\theta_1+\phi+\alpha)\Gamma(2\alpha)}
 {\Gamma(\theta_1+\phi+2\alpha)\Gamma(\alpha)} +
\frac{\Gamma(\theta_2+\phi+\alpha)\Gamma(2\alpha)}{\Gamma(\theta_2+\phi+2\alpha)\Gamma(\alpha)}. \label{eq:mbet}
\end{equation}
Determining the function $\alpha(\bftheta)$ explicitly is therefore not possible in general.
We can make progress when $\alpha=2$. In this case \eqref{eq:mbet} becomes
\[ m(\bftheta,\phi) = \frac{6}{(\theta_1+\phi+3)(\theta_1+\phi+2)} + \frac{6}{(\theta_2+\phi+3)
(\theta_2+\phi+2)}, \]
and from this
\[ \alpha(\bftheta) = -\frac12(\theta_1+\theta_2) - \frac52 + \frac12 \sqrt{(\theta_1-\theta_2)^2+
2\sqrt{25(\theta_1-\theta_2)^2+144} +25}. \]
Minimizing $\bfa.\bftheta+\alpha(\bftheta)$ by setting $\phi=\theta_1+\theta_2$ and $\psi=\theta_1-\theta_2$ 
we see that the only values which are greater than $-\infty$ are on the line $a_1+a_2=1$ and they are found as
\[  - \frac52 +\inf_{\psi}\left( \frac{(a_1-a_2)\psi}{2} +  \frac12 \sqrt{\psi^2+ 2\sqrt{25\psi^2+144} +25}
\right).\]
This does not have an explicit solution but can be computed numerically.

Numerical solutions could also be found for other integer values of $\alpha$, as the expressions 
in \eqref{eq:mbet} become polynomial.

In the limit as $\alpha \to \infty$ we have a deterministic split in that the rectangle is always divided 
in half and only the direction of the split is random. In this case we can compute $\alpha^*$. Note that we have
$E(U^{\gamma})=2^{-\gamma}$ which gives
\[ \alpha(\bftheta) = \frac{\ln(2^{-\theta_1}+2^{-\theta_2})}{\ln{2}}-1. \]
Thus we can compute $\alpha^*(\bfa)$ as
\begin{eqnarray}\label{1809141037}
\alpha^*(\bfa) = \left\{ \begin{array}{ll} \displaystyle{\frac{a_1\ln{a_1} + a_2\ln{a_2}}{\ln{(1/2)}}} & a_1+a_2=1 \\ 
-\infty & a_1 +a_2 \neq 1. \end{array} \right. 
\end{eqnarray}
As we have a growth function which exists only on a line segment in $\mathbb{R}_+^2$, we can give results on the 
exponential growth of the number of rectangles of a certain size in the model using Corollary~\ref{cor:gbrw}.

It is not difficult to see that we have a continuous dependence on $\alpha$ in this model, even though 
analytically tractable $\alpha$ are rare! Thus we could envisage using this as a basis for a statistical 
estimation procedure to determine the parameter $\alpha$ from data on the number of rectangles of given 
length and height.

(P) In this setting we can find the power law exponent in the biased case. For the case of $\alpha=2$ 
we see that the Malthusian parameter in the branching process will satisfy
\begin{equation}
 1= E \sum_{i=1}^2 U_i^{\gamma} \frac{T_0^i+1}{2}. \label{eq:gamparam}
 \end{equation}
We can use the independence of $U_i$ and $T_0^i$ to see that
\[ 1 = \frac{(ET_0+1)6}{(\gamma+3)(\gamma+2)}, \]
and hence as $ET_0=1/(1-2p)$,
\[ \gamma = -\frac52+\frac12 \sqrt{\frac{49-50p}{1-2p}}. \]
In the general case, to find the exponent $\gamma$, we need to solve 
\begin{equation}
1 = \frac{2-2p}{1-2p}\frac{\Gamma(\gamma+\alpha)\Gamma(2\alpha)}{\Gamma(\gamma+2\alpha)\Gamma(\alpha)}. 
 \label{eq:gam}
\end{equation}
By continuity of the Gamma function the exponent $\gamma$ will be a continuous function of $p,\alpha$ for 
$0<p<1/2$ and $\alpha>0$. Also in the case where $\alpha\to\infty$, using $E(U^{\gamma})=2^{-\gamma}$, 
we can solve equation 
\eqref{eq:gamparam} explicitly to get
\begin{eqnarray}\label{1809131723}
 \gamma = \frac{\ln{\frac{2-2p}{1-2p}}}{\ln{2}}. 
\end{eqnarray}
A graph of the value of $\gamma$ as a function of $\alpha$ for two different values of $p$ is shown 
in 
~\ref{1809131045}-Right.


As $p\to 1/2$ we can observe the asymptotics in the exponent $\gamma(p)$ as a function of the parameter $\alpha$. 
It is clear from equation \eqref{eq:gam}
that we want $\gamma(p)$ such that
\[ \frac{\Gamma(\gamma(p)+2\alpha)}{\Gamma(\gamma(p)+\alpha)} = \frac{2-2p}{1-2p} 
\frac{\Gamma(2\alpha)}{\Gamma(\alpha)}. \]
Note that the right hand side of this equation diverges as $p\to 1/2$ and hence we require the large parameter 
asymptotics of the Gamma function;
\[ \lim_{x\to\infty} \frac{\Gamma(x+a)}{\Gamma(x)x^a}=1. \]
Employing this we see that
\[ \lim_{p\to 1/2} \gamma(p)^{\alpha}(1-2p) = \frac{\Gamma(2\alpha)}{\Gamma(\alpha)}=C_{\alpha}, \]
and hence
\[ \gamma(p)\sim \left(\frac{C_{\alpha}}{1-2p}\right)^{1/\alpha}\;\;\mbox{as }p\to 1/2. \]
In the case where $\alpha=\infty$, from the explicit expression, we see that $\gamma(p)$ diverges 
logarithmically as $p\to 1/2$.

%
%
%

\section{Numerical results}\label{1801071851}

We present a set of numerical results for realizations of
the fragmentation schemes discussed throughout this paper and 
compute and analyze the statistics of the patterns generated.
All the numerical routines are implemented in Matlab and make use of Matlab's pseudorandom number generator. 

\subsection{Asymptotics of rectangles}\label{1708281611}

\begin{figure}[h!]
    \centering
    \begin{subfigure}[b]{0.31\textwidth}
\includegraphics[width=6.2cm]{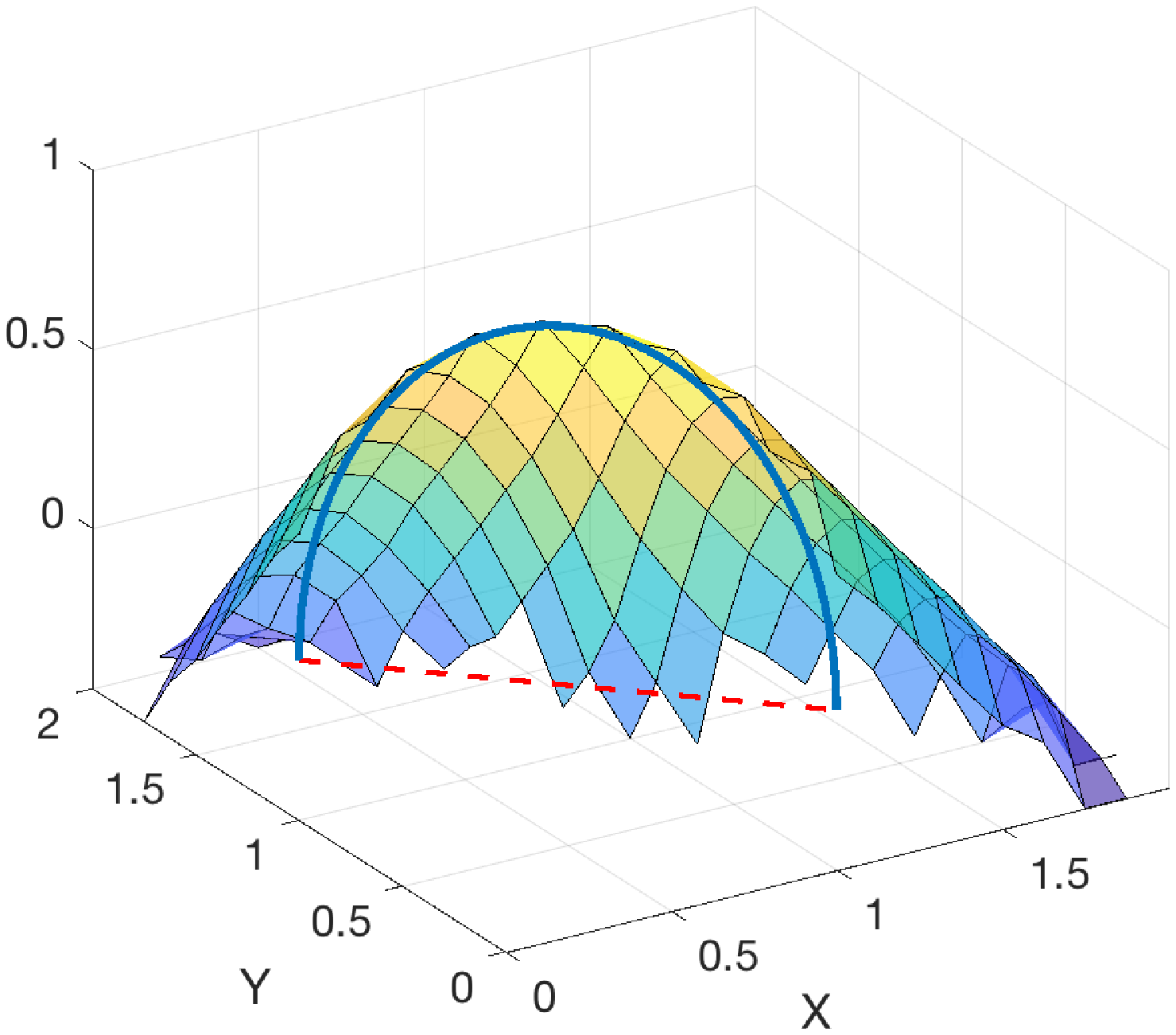}
        \caption{}\label{p05}
    \end{subfigure}
\hspace{-0.0cm}
    ~ 
    \begin{subfigure}[b]{0.31\textwidth}
\includegraphics[width=6.2cm]{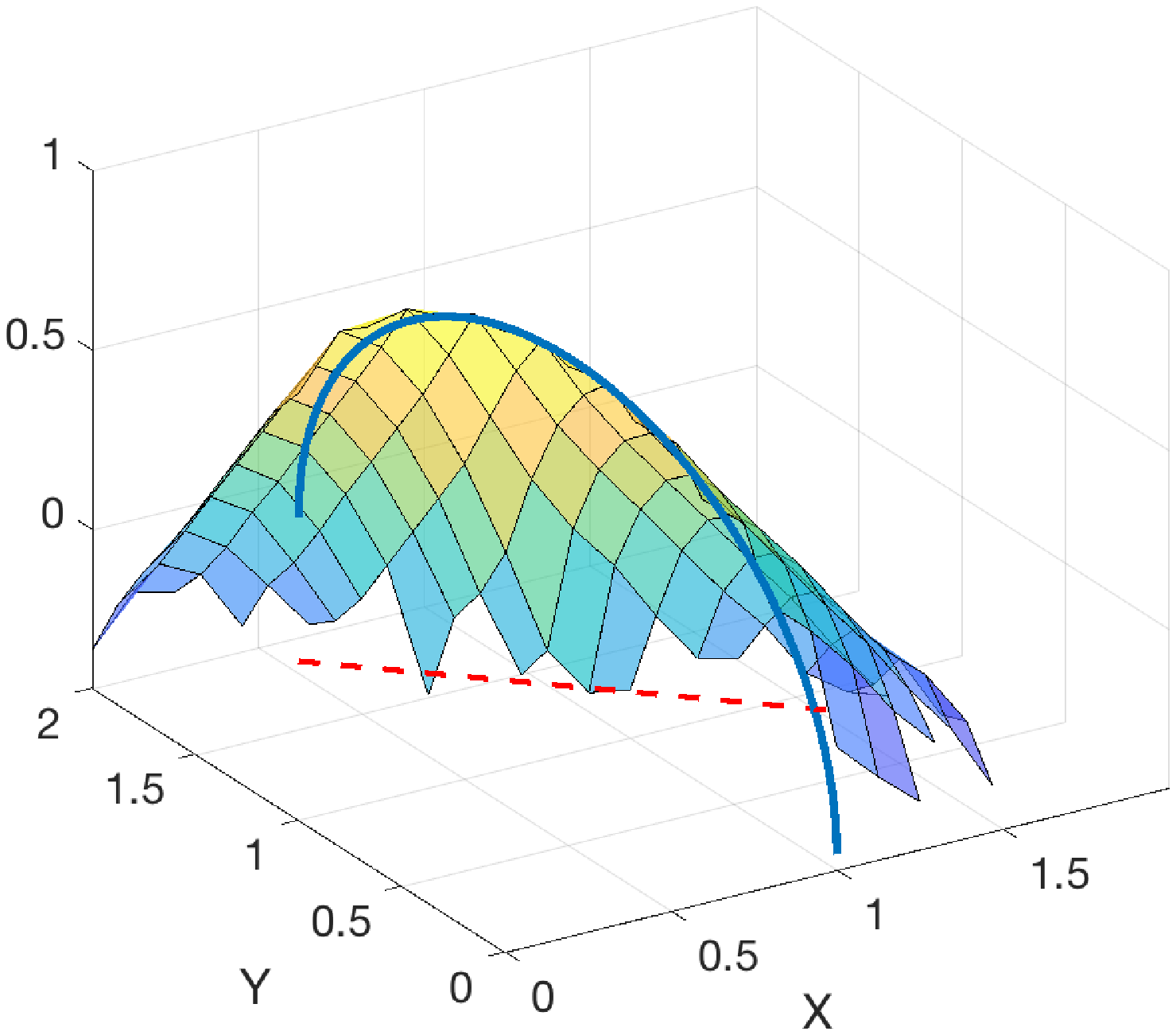}
        \caption{}
        \label{p07}
    \end{subfigure}
\hspace{-0.0cm}
    \begin{subfigure}[b]{0.31\textwidth}
\includegraphics[width=6.2cm]{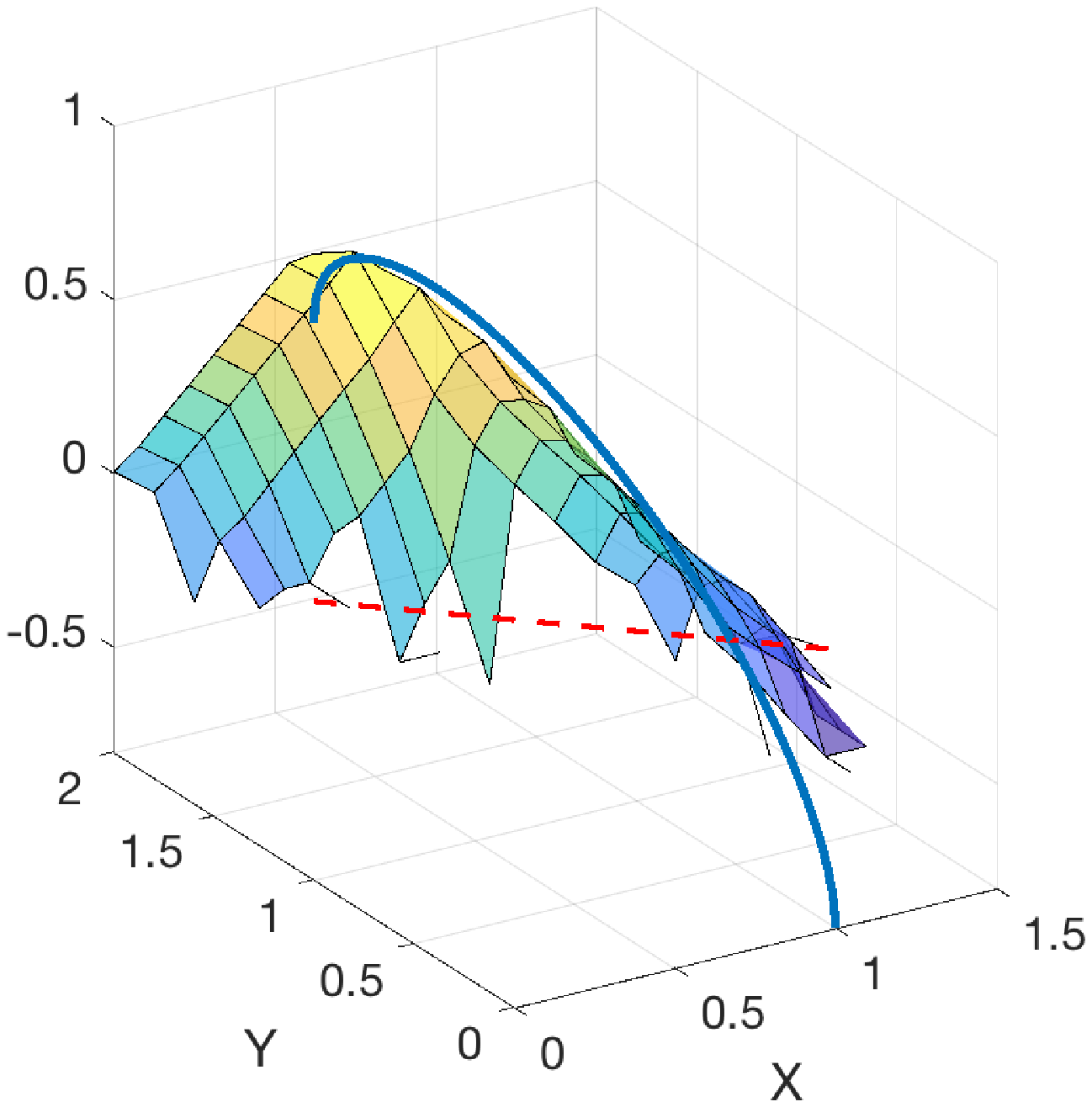}
        \caption{}
        \label{p09}
    \end{subfigure}
\caption{
Averaged histograms for the distribution of rectangles obtained for a fragmentation model parametrized by different values of $p$. 
From left to right: $p=0.5,0.3$ or $0.1$ respectively. The exact analytical solution $f_p$ is represented by a continuous curve.
The set $x+y=1$ corresponds to the dashed line.
}\label{cartoon1}
\end{figure}

Theorem \ref{1708261300} provides us with the asymptotic description of the rectangles generated in the 2D 
fragmentation process described in Section~\ref{sec:model}. Although this works in the limit as $t\to\infty$, we now show how to compare the analytical result with numerical computations of the distribution of rectangles obtained for finite $t$.
The average distribution of rectangles will be represented in the form of a 3-dimensional histogram computed 
by binning the rectangles based on the length of their sides. Each rectangle of coordinates $(a,b)$ with 
$0<a,b<1$ contained in the unit square is mapped onto the space of normalized coordinates $(x,y)$ defined in Theorem \ref{1708261300}.
%
%
To fix ideas, consider the unbiased case with $p=1/2$. In Figure~\ref{p05} we report the averaged logarithmic
histogram of the distribution of the rectangles stopped at $n=2\times 10^5$ which corresponds to 
$t\approx 11.51$. Indeed, $t$ grows logarithmically as a function of $n$.
An averaged histogram is obtained by averaging out 100 realizations of the 2D fragmentation
process and it is represented in a logarithmic space. For finite time, the distribution corresponds to a surface 
in the space $(x,y)$ and, in the limit as $t$ tends to infinity, the surface collapses onto a curve defined 
over the line of equation $x+y=1$.
%
%
For general $p$, the exact formula $\alpha^* $ for the asymptotic average distribution of the rectangles 
is given in (\ref{1708281411}).
%
In order to make the dependence on the parameter $p$ explicit and make the notation more intuitive, 
here we identify $x\equiv a_1,y\equiv a_2$  and define the new function $f_p(x,y):=\alpha^*(a_1,a_2)$ 
which, in turn, yields 
$$
f_p(x,1-x)= 4\sqrt{p(1-p)}\sqrt{x(1-x)}
+(2p-1)(2x-1),\qquad 0\leq x\leq 1,
$$
when $x+y=1$, while we assume $f_p\equiv-\infty$ if $y\neq 1-x$.
%
%
Observe that the shape of the histogram 
(and, in particular, the maximum point) depends heavily on a set of parameters such as the number, size 
and location of the bins.
Here, we decide to normalize the averaged histogram so that its maximum value coincides with that of the 
analytical solution, even for a finite $t$. By doing this  we implicitly assume the only relevant  information 
is the shape of the histogram, which can be directly compared with the profile of the analytical solution.
Comparison with the analytical solution shows good matching of the profiles in all the cases under 
consideration \textcolor{black}{~\ref{p05}-\ref{p09} }.
From an inspection of the histograms
we note that at finite time the majority of the rectangles tend to concentrate along the line of equation
$x+y= 1$. In the untransformed variables, this is equivalent to the curve given by the equation $ab = e^{-t}$, 
the set of the rectangles of areas equal to  $ e^{-t}$, which is in agreement with the assumptions of a model 
for which the individuals alive correspond to rectangles with area less than $ e^{-t}$ at time $t$.


\subsection{Power laws}\label{1708261900}
 
Focusing on the alternative viewpoint for counting interfaces described in Section \ref{1708081205}, we consider the statistics of interfaces 
leading to power laws. Beginning with the 2-dimensional case, we denote by $N^p_h(x) $ the cumulative distribution of  horizontal interfaces of 
length greater or equal than $x$ in a biased fragmentation process, where  $0<p<1$ is the probability of an interface growing horizontally.
%
%
From Theorem \ref{1708280114} we deduce $N^p_h(x)$ has a power law profile with exponent equal to $-(\frac{1}{1-2p})$ when $p<1/2$. 
Analogously, the  derivative of $N^p_h(x) $, regarded as the \textit{density} of the horizontal interfaces of length $x$ (in other words, the 
number of  interfaces of length contained in the interval $(x,x+dx)$) displays power law behavior  with characteristic exponent equal to 
$-1-(\frac{1}{1-2p})$.

\begin{figure}[h!]
\centering%
\includegraphics[width=7cm]{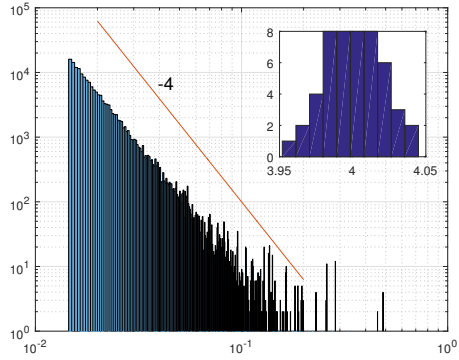}
\caption{
Histogram of the density of the areas of horizontal plates in a 3D fragmentation model.
For this particular realization we have counted $166866\approx 5\times 10^5\times \frac{1}{3}$ horizontal interfaces. 
A line of slope $-4$ has been drawn in the log-log plane to guide the reader's eye. 
We plot (in-picture) the histogram for the (absolute values of the) exponents of the densities of 
interfaces obtained in 50 realizations of a 3D fragmentation scheme, each of which is composed of 
$5\times 10^5$ cuboids. Computed average and standard deviation are reported in Table \ref{1708280154}.
}\label{1708171320}
\end{figure}

\begin{figure}[h!]
\centering%
\includegraphics[width=3cm]{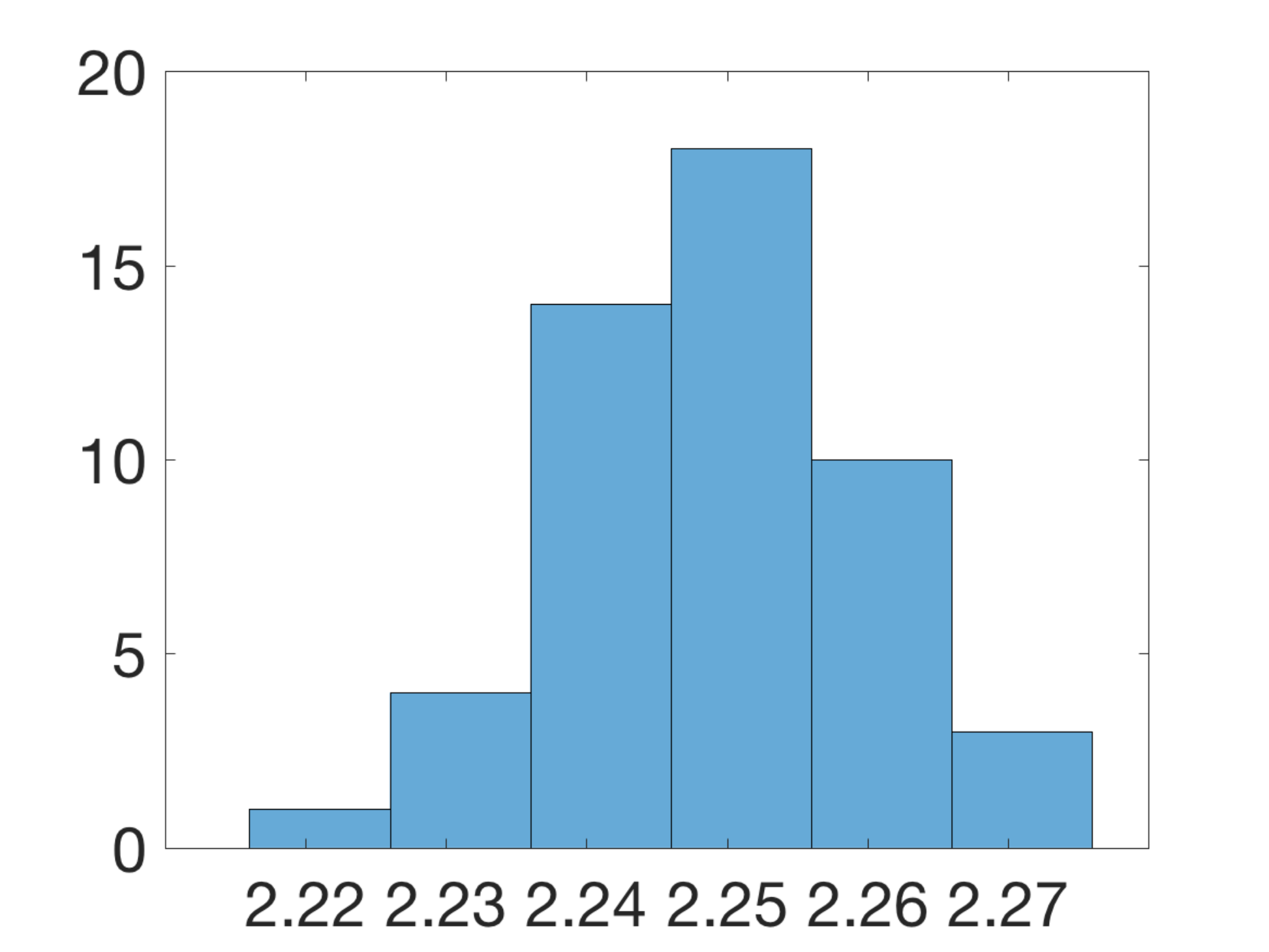}
\includegraphics[width=3cm]{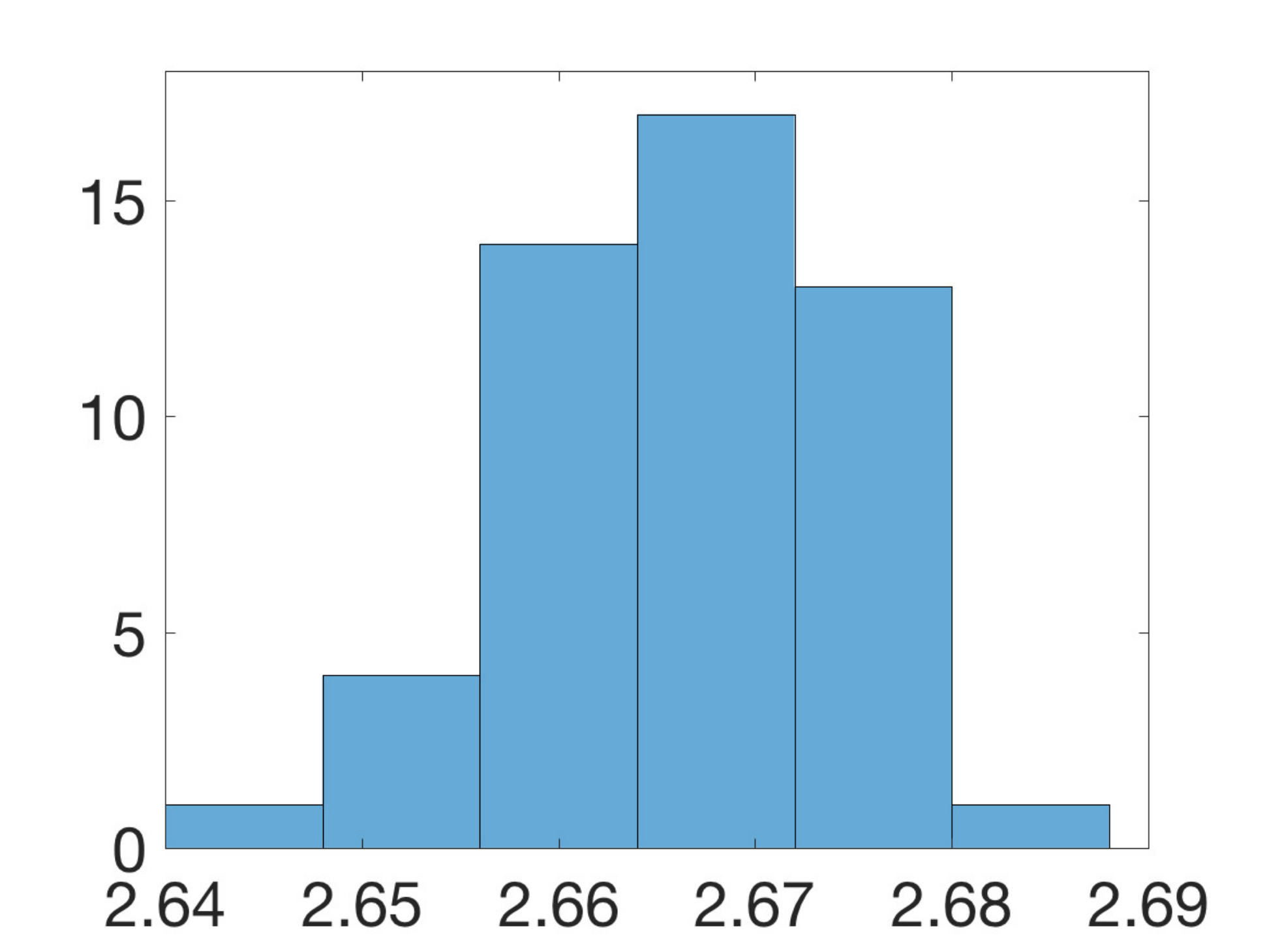}
\includegraphics[width=3cm]{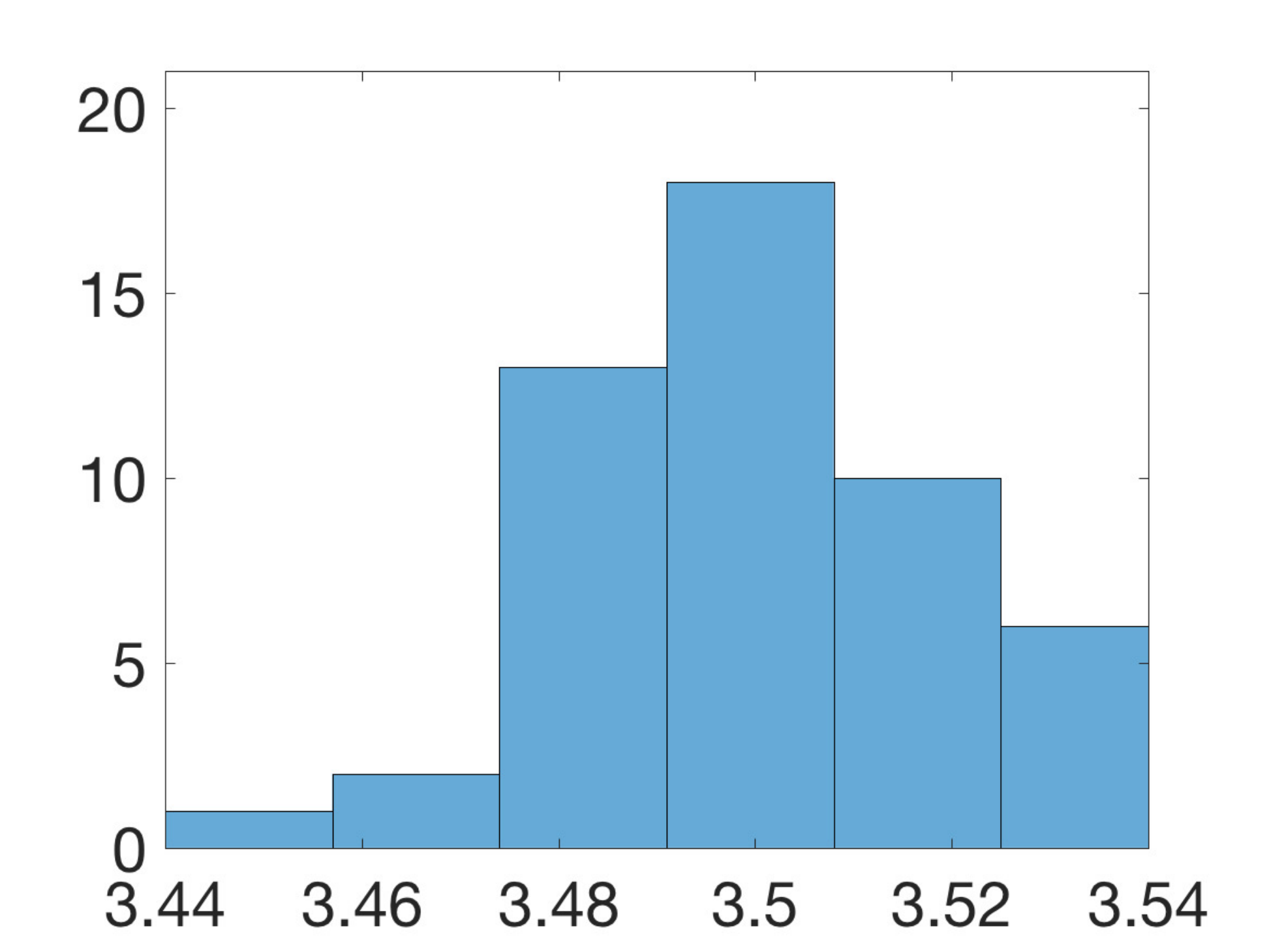}
\includegraphics[width=3cm]{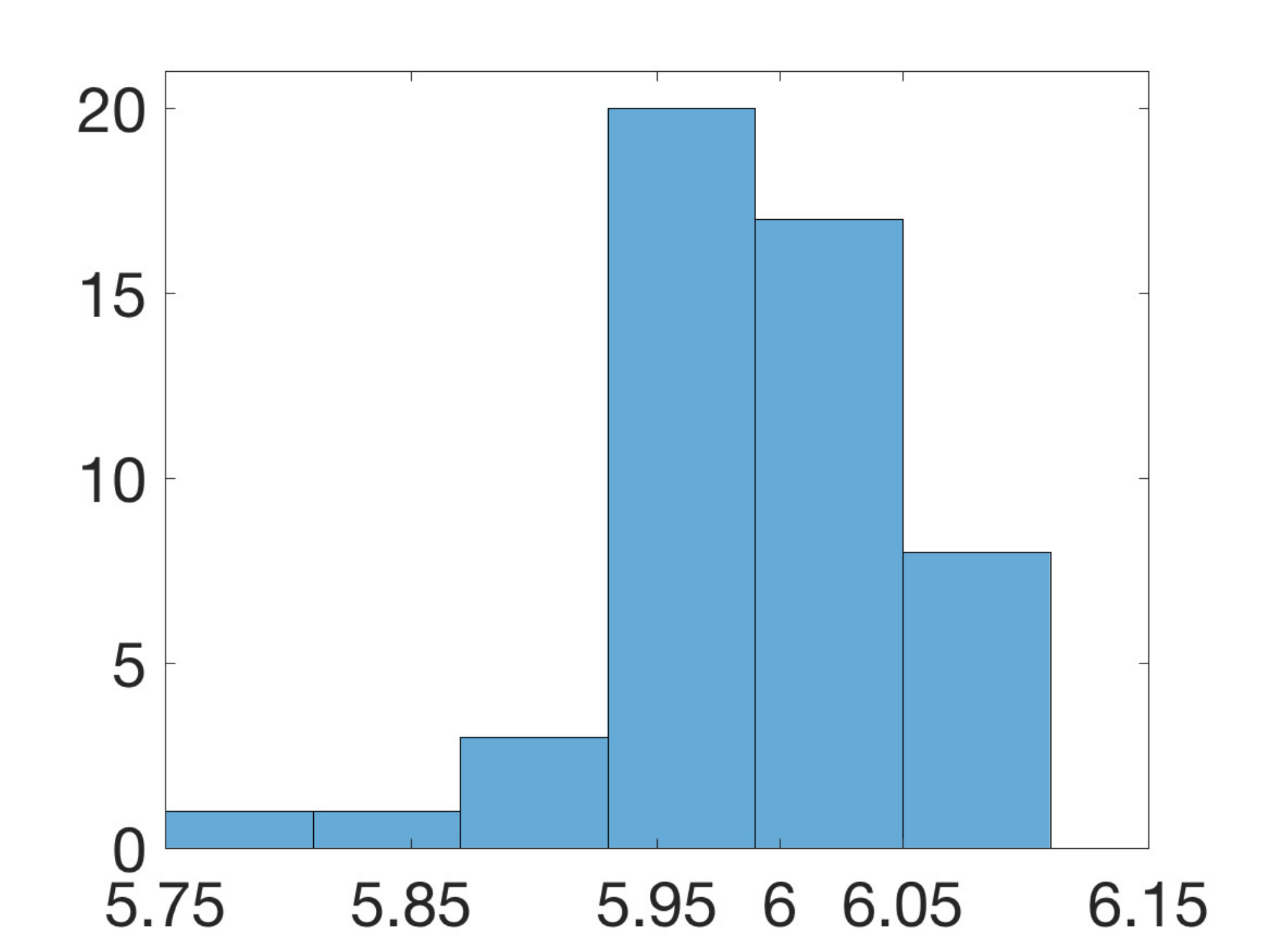}
\includegraphics[width=3cm]{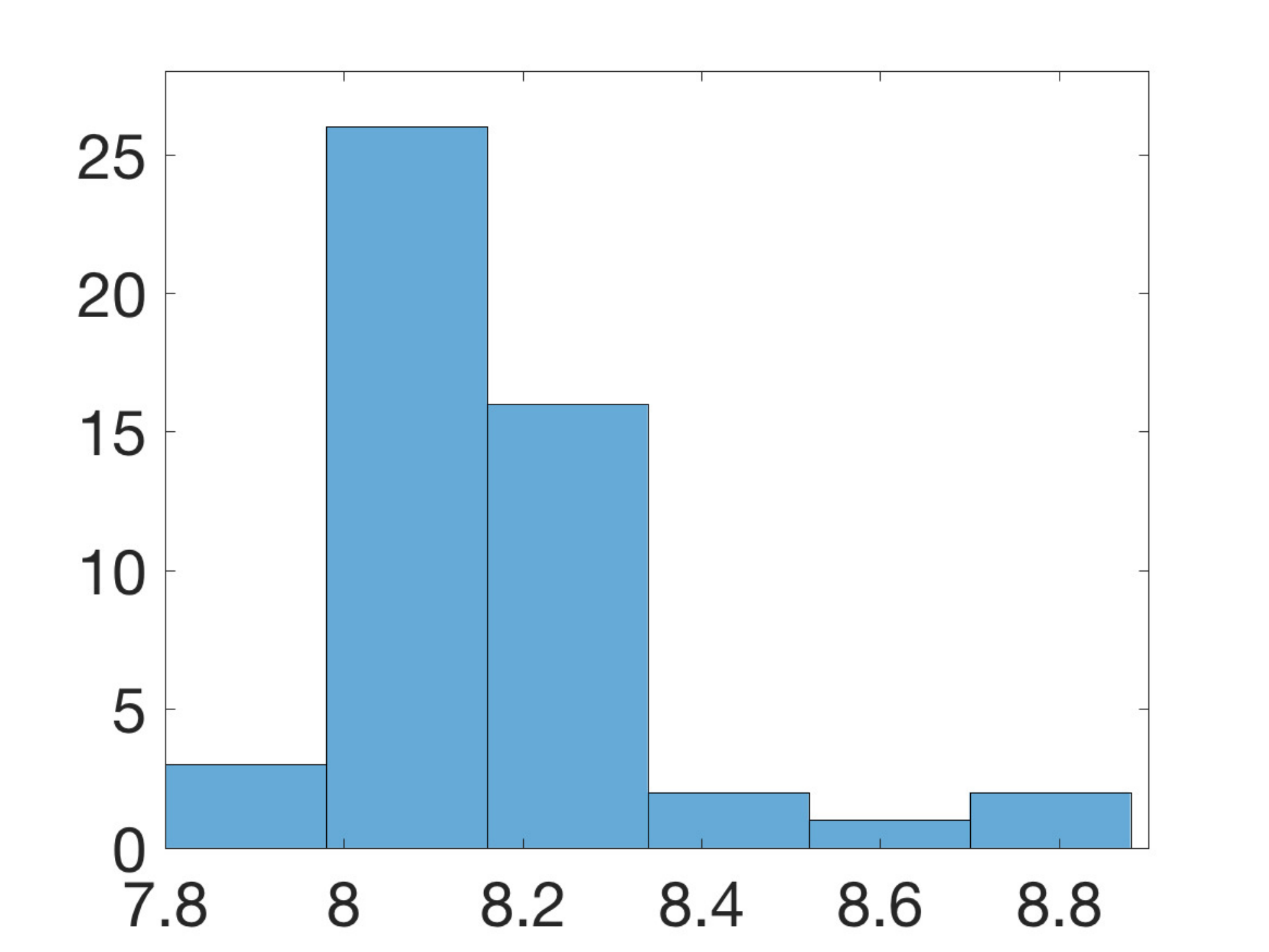}
\caption{
Histograms  of the (absolute values of the) exponents of the densities associated with the 2D fragmentation 
process parameterised by $p$. 
Left to right: the parameter $p$ ranges through $\{0.1,0.2,0.3,0.4,0.43\}$  with the corresponding 
predicted exponents equal to $\{-2.25, -2.\overline{6}, -3.5, -6, -8.142...\}$ respectively.
Each  histogram has been constructed from 50 realizations of the 2D version of the fragmentation process 
with each realization containing $3\times 10^5$ rectangles. 
Computed averages and standard deviations are reported in Table \ref{1708280154}.
}
\label{1708212152}
\end{figure}

The analytical prediction for the exponent has been tested with simulations of fragmentation schemes analogous to those described in 
Section~\ref{sec:model}. The procedure consists of assigning a value to $p$ and then running several realizations of the fragmentation routine, 
each one composed of $n=3\times 10^5$ interfaces for each value of $p$. To make our routine faster, we split at each iteration the rectangle 
with largest horizontal side first so that at each step a significant amount of data is available for our analysis.
The numerical results confirm the histograms of the density of the interfaces obtained in this way do have a power law behaviour (see, for instance, 
Figure \ref{1708171320}). In order to determine the exponents of the power laws, we have employed Clauset's library \cite{Clauset09}
based on the Maximum Likelihood Method (MLM). We remark the MLM is insensitive to the histogram parameters, such as the number, 
size and location of the bins and is stable with respect to possible fluctuations in the tail of the distribution.
The results reported in Table  \ref{1708280154} and in Figure \ref{1708212152} show in general excellent agreement with those 
predicted by the theory, even for $p$ close to 1/2. As $p\to 1/2$ the number of horizontal interfaces of length less than or equal to $x$ tends to 
infinity, and in turn, the exponent of the power law for the density of horizontal interfaces diverges. 
We pushed our numerical experiments up to  
case $p=0.43$, in which the density of horizontal interfaces is characterized by an exponent which is approximately  $-8$. Even 
in this situation we obtain an overall good agreement of the computed averaged exponent although with much 
less accuracy than for smaller $p$'s.

In the 3D version of the fragmentation process, a unit cube is fragmented by planar interfaces nucleating and growing perpendicular to the 
three coordinate axes with equal probability, that is,  with $p_1=p_2=1/3$ as described in Section~\ref{sec:model}. If we focus on the 
horizontal plates (or, in other words, the ones perpendicular to the direction $e_3$ in the Cartesian frame), we expect the exponent for the density of the plates of area equal to $x$ to coincide with the exponent for the density of horizontal lines equal to $x$ in a 
biased 2D process with $p=1/3$, which is $-1-(\frac{1}{1-2/3})=-4$. In  Figure \ref{1708171320} we show the histogram of the density
of the areas of the horizontal interfaces in one realization of a 3-Dimensional fragmentation process. The histogram shows a clear power law 
behaviour with a sharp $-4$ exponent, in agreement with the prediction.
In order to compute the statistics of the exponents of the density of the areas of the horizontal plates, a series of 50 realizations of 
the 3D fragmentation were run, with each realization containing $5\times 10^5$ plates. The averaged exponent and its standard deviation from the 50 realizations are displayed in Table \ref{1708280154}, showing excellent agreement with the prediction.

\begin{table}[h]
\centering
    \begin{tabular}{| l | l | l | l | l | l |l | }
    \hline
    $p$ & 0.1 & 0.2 &0.3 &1/3$^{\star}$ & 0.4 &0.43 \\ \hline\hline
    Predicted exponent & -2.25   &  -8/3   & -3.5 & -4  & -6  & -8.142..\\ \hline
    Computed exponent & -2.2493    & -2.6664    &  -3.5012 & -4.0013  &   -5.9894    &  -8.1624 \\ \hline 
 Standard deviation  & 0.0106  &   0.0079    & 0.0194 &  0.0193 &  0.0604    &      0.1932 \\ \hline
    \end{tabular}
\caption{Analytical and numerical computation of the power law
exponents for the densities of the length of horizontal interfaces in a biased fragmentation process parameterized by $p$.
The computed exponents are averages over 50 realizations of a fragmentation process, each of which is composed of $3\times 10^5$ rectangles.
$^{\star}$The case corresponding to $p=1/3$ refers to the 3D model. Here the data correspond to the density of the areas of horizontal plates in an unbiased fragmentation process. Averages are computed over 50 realizations of a fragmentation process, each of which is composed of $5\times 10^5$ cuboids.
}
\label{1708280154}
\end{table}
 
%
%
%
%

\subsection{Constant rate splitting}

We focus on the approach described in Section~\ref{sec:indep} consisting of splitting rectangles without 
regard to their size. The fragmentation mechanisms yield highly inhomogenous microstructures and consequently
distribution profiles for the rectangles 
which are very spread out,
a striking difference from the case analysed 
in Section \ref{sec:largest}, where at each step the largest rectangle splits.
%
%
%
%
%
%
%
%
%

\begin{figure}[h!]
\centering%
\includegraphics[width=7cm]{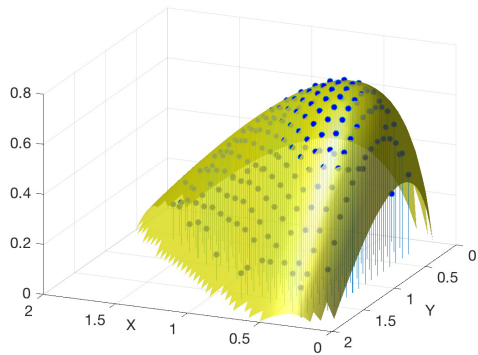}
\caption{Uniform splitting fragmentation process.
A normalized histogram for the distribution of rectangles of 
one realization of the fragmentation mechanism 
occurring via discrete generations and
described in Section \ref{201809141440}. The histogram is represented here by a stem-plot. 
 Binning is uniform in the transformed 
space with bin size equal to 0.1 in both directions. 
The microstructure obtained is composed of 
$2^{21}$ rectangles. The full-colour surface 
corresponds to the analytical solution.
}
\label{1706211243}
\end{figure}

Consider the
%
discrete-generation mechanism
of Section \ref{201809141440}.
%
An example of a microstructure obtained according to this method is displayed in Fig. \ref{2006172133}-Right.
Equation (\ref{1708262241}) provides us with the expected asymptotic distribution 
of the rectangles obtained according to this fragmentation model. 
The general rectangle with coordinates $(a,b)$ is mapped into the space of normalized logarithmic 
coordinates $(x,y)$ and  time $t=\log_2 (2^n)=n$ coincides with the number of generations.
The exact formula for the asymptotic density of the rectangles is then given by
\begin{eqnarray}\label{1803251337}
g(x,y)=1-x-y+2\log(\sqrt{x}+\sqrt{y})\qquad x,y>0.
\end{eqnarray}
In particular, observe that $g$ attains its maximum value at the point $(x,y)=(1/2,1/2)$, where 
$g(\frac{1}{2},\frac{1}{2})=\log 2$.
When running a set of realizations of this fragmentation scheme for subsequent generations (in other words, for increasing values of $n$) and if we construct the corresponding histograms
in the logarithmic space $x,y$ while keeping
the number, location and size of the bins
constant we observe
the profile of the histogram rises and, at the same time, the coordinates of the maximum point 
move towards $(1/2,1/2)$.
%
%
%
%
We turn this observation into a method for 
a more effective comparison of the 
histogram and the analytical solution.
%
%
%
%
%
%
%
%
%
%
%
%
The profile of the logarithmic histogram for the distribution of the rectangles obtained in a single realization 
of the fragmentation process is shown in Figure \ref{1706211243}. It has been normalized and shifted rigidly
so that both its maximum and the location of its maximum 
correspond to those of $g$. We are therefore left with comparing the shape of the analytical solution 
of the asymptotic density of the rectangles
and the profile of the histogram for a single (finite) microstructure which show a reasonably good agreement.

\subsection{Beta distribution}

\begin{figure}[h!]
\centering%
\includegraphics[width=7cm]{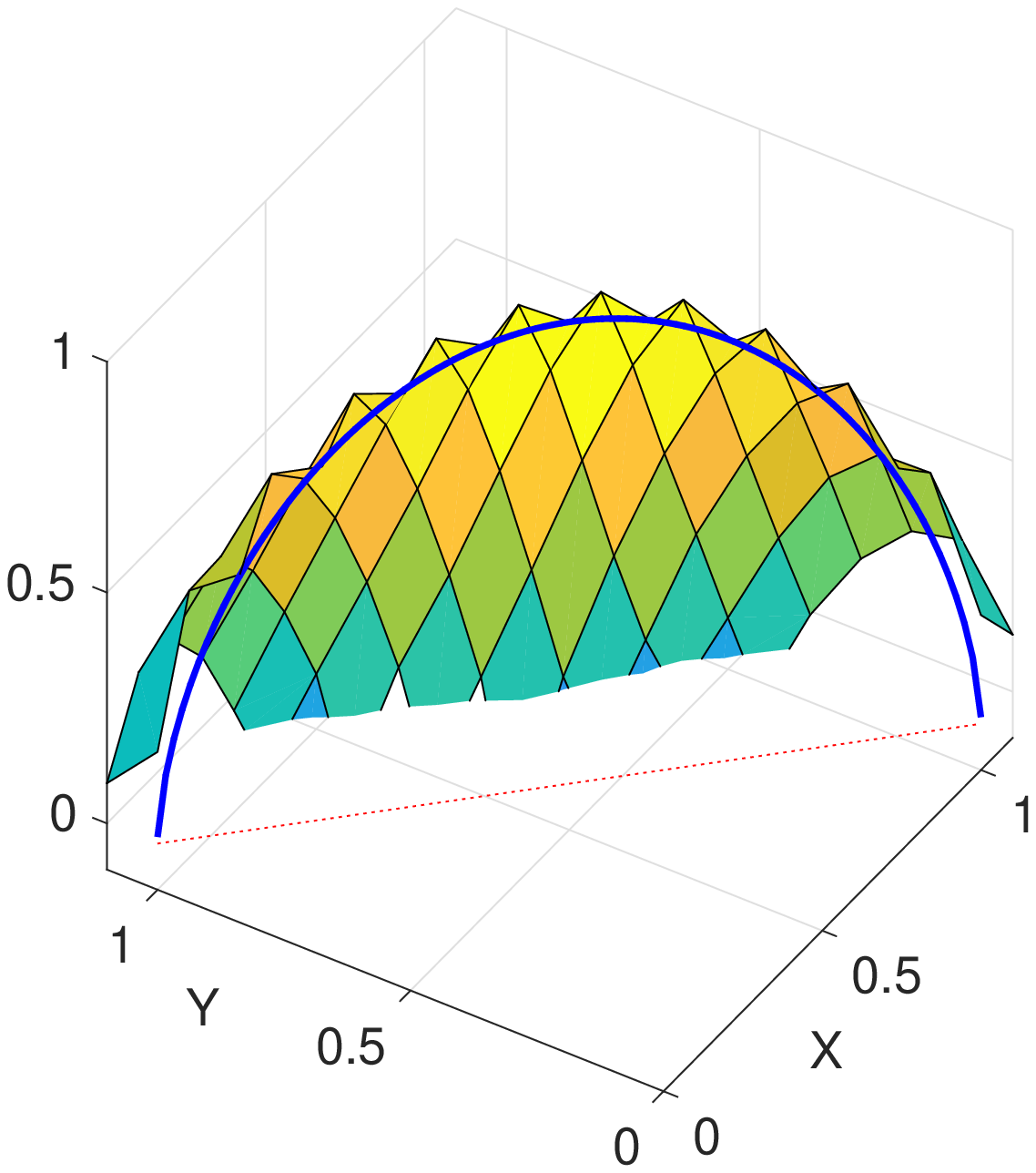}
\hspace{-1.0cm}
\includegraphics[width=9cm]{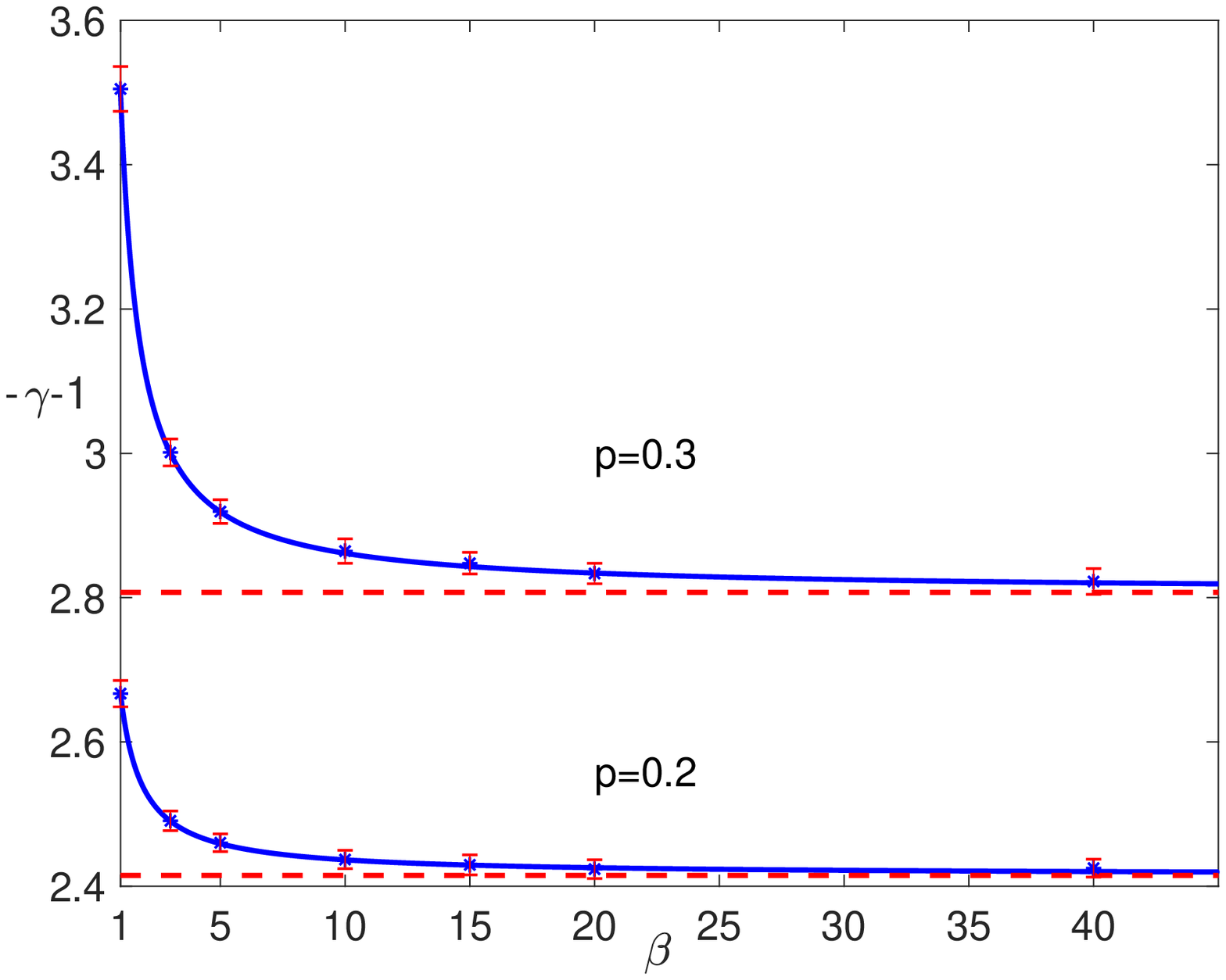}
\caption{Fragmentation experiments where nucleation follows a beta distribution process as described in Section \ref{1809241600}.
Left: Surface plot of the normalized and averaged histogram of the distribution of rectangles for an unbiased process ($p=\frac{1}{2}$). The analytical solution is represented by a continuous curve.
The set $x+y=1$ corresponds to the dotted line.
Right: Power law exponents of the density of horizontal interfaces obtained for a biased process with $p=0.3$ or $0.2$, respectively. Continuous curves represent the value of the exponents
as a function of $\beta$.
Dashed lines represent the asymptotic exponents 
obtained when $\beta\to\infty$.
The average exponents are
obtained by running $100$ realizations of the fragmentation process (each of which is composed of $10^5$ rectangles) and error bars (measuring standard deviation) are displayed.
}
\label{1809131045}
\end{figure}

Numerical computations for the distribution of rectangles obtained when the nucleation occurs according to a beta distribution show excellent agreement with analysis.
A plot for the averaged histograms is displayed in Figure \ref{1809131045}-Left and compared with the analytical solution available for $p=\frac{1}{2}$ given by Eq. (\ref{1809141037}) (approach (G) in Section \ref{1809241600}).
Here 100 realizations of the fragmentation have been performed, each containing $2\times 10^4$ rectangles.
For details on how the normalized histogram has been constructed see Section \ref{1708281611}.

Following approach (P) in Section \ref{1809241600}, power laws exponents for the density of the length of horizontal interfaces have been computed and compared with the analytical prediction (Figure \ref{1809131045}-Right). 
Recall $-\gamma$ is the exponent for the cumulative distribution of the horizontal interfaces. 
The continuous lines in Figure (\ref{1809131045}-Right)
correspond to the values of $-\gamma-1$, where $\gamma$ is obtained by inverting
(\ref{eq:gam}) as
a function of $\beta$ and where $p$ is regarded as a parameter (equal to either $0.3$ or $0.2$).
Numerical results for the average exponents (along with their corresponding standard deviation) are obtained from a set of $100$ realizations of the microstructure, each containing $10^5$ rectangles. 
The limit cases correspond to $\beta\to\infty$ (for which 
$\gamma$ can be computed explicitly from (\ref{1809131723})) and $\beta=1$ (in which case the beta distribution coincides with the uniform distribution and therefore $\gamma\equiv\frac{1}{1-2p}$, see Theorem~\ref{1708280114}).

\paragraph{Acknowledgments.}

We thank Prof. John Ball for suggesting the model and Prof. Eduard Vives for discussions about 
the TPV model.
PC is grateful to Prof. Giovanni Zanzotto for insightful discussions.
We also thank Dr. Robert Niemann for kindly giving us permission to publish 
the illustrative image in Figure~\ref{1609181210}.
PC is supported by JSPS Research Category	
Grant-in-Aid for Young Scientists (B) 16K21213 and an AiRIMaQ Fellowship. 
PC is grateful to the hospitality of Mathematical Institute, University of Oxford. 
PC holds an honorary appointment at La Trobe University and is a member of GNAMPA.
 

%


\begin{thebibliography}{10}



\bibitem{JB02}
Ball, J.M.  Some open problems in elasticity. In Geometry, Mechanics, and Dynamics, pages 3--59, Springer, 
New York, 2002.

\bibitem{BCH}
Ball, J.M., Cesana, P.and Hambly, B.M. \emph{A probabilistic model for martensitic avalanches}. 
MATEC Web of Conferences Volume: 33 Article Number: 02008, 2015.

\bibitem{Ball87}
Ball, J.M. and James, R.D. \textit{Fine phase mixtures as minimizers of energy}, Arch. Rat. Mech. Anal., 
{\bf 100} (1987), 13--52.

\bibitem{Ber}
Bertoin, J. \emph{Random fragmentation and coagulation processes} Cambridge Studies in Advanced Mathematics, 
102. Cambridge University Press, Cambridge, 2006.

\bibitem{Bhatta}
Bhattacharya, K. {\it Microstructure of Martensite}, Oxford University Press, 2003.

\bibitem{Big78}
Biggins, J. D. The asymptotic shape of the branching random walk. \emph{Advances in Appl. Probability}
  {\bf 10} (1978), 62--84. 
  
\bibitem{Big90}
Biggins, J. D. Uniform convergence of martingales in the branching random walk. \emph{Ann. Probab.} 
{\bf 20} (1992), 137--151.

\bibitem{Big94}
Biggins, J. D. How fast does a general branching random walk spread? \emph{Classical and modern branching 
 processes (Minneapolis, MN, 1994)}, 19--39, IMA Vol. Math. Appl., 84, Springer, New York, 1997.

\bibitem{Big95}
Biggins, J.D. The growth and spread of the general branching random walk,
\emph{Ann. Appl. Probab.}, {\bf 5} (1995), 1008--1024.

\bibitem{BNS}
Broutin, N., Neininger, R. and Sulzbach, H. Partial match queries in random quadtrees. \emph{Proceedings of the 
Twenty-Third Annual ACM-SIAM Symposium on Discrete Algorithms}, 1056–1065, ACM, New York, 2012.

\bibitem{CCH} Charmoy, P.H.A, Croydon, D.A. and Hambly, B.M.  Central limit theorems for the spectra of a class of random self-similar fractals. Trans. Amer. Math. Soc.{\bf 369} (2017), 8967--9013. 

\bibitem{Clauset09}
Clauset, A., Shalizi, C.R. and Newman, M.E.J. Power-Law Distributions in Empirical Data
\emph{SIAM Rev.}, {\bf 51} (2009), 661--703. 
Additional material available online at \url{http://tuvalu.santafe.edu/~aaronc/powerlaws}

\bibitem{Derrida87}
Derrida, B. and Flyvbjerg, H. Statistical properties of randomly broken objects and of multivalley structures 
in disordered systems \emph{J. Phys. A}, {\bf 20} (1987), 5273--5288.

\bibitem{Frontera95}
Frontera, C., Goicoechea, J., R\'afols, I. and Vives, E. Sequential partitioning: An alternative to understanding 
size distributions of avalanches in first-order phase transitions. \emph{Phys. Rev. E} 52, 5671 (1995)

\bibitem{GST}
Georgii, H.-O.; Schreiber, T. and Thäle, C. Branching random tessellations with interaction: a thermodynamic view. 
\emph{Ann. Probab.} {\bf 43} (2015), 1892--1943.

\bibitem{HKK}
Harris, S. C., Knobloch, R. and Kyprianou, A. E. Strong law of large numbers for fragmentation processes. 
\emph{Ann. Inst. Henri Poincaré Probab. Stat.} {\bf 46}  (2010),  119--134.

\bibitem{MM} 
Mackisack, M.S. and Miles, R.E. Homogeneous random tessellations. \emph{Adv. Appl. Probab.} {\bf 28} 
(1996), 993--1013.

\bibitem{Muller}
M\"{u}ller, S. ``\emph{Variational methods for microstructure and phase transitions}'', in: 
Proc. C.I.M.E. summer school ``Calculus of variation and geometric evolution problems'', 
Cetraro 1996, (F. Bethuel, G. Huisken, S. M\"{u}ller, K. Steffen, S. Hildebrandt, M. Str\"{u}we eds.),
Springer LNM vol. 1713, 1999.

\bibitem{NW}
Nagel, W. and Weiss, V. Crack STIT tessellations: characterization of stationary random tessellations 
stable with respect to iteration. \emph{Adv. in Appl. Probab.} {\bf 37} (2005), 859--883.

\bibitem{Nandi16}
Nandi, S.K., Biroli, G. and Tarjus, Gi. Spinodals with Disorder: From Avalanches in Random Magnets to Glassy 
Dynamics \emph{Phys. Res. Let.} {\bf 116} (2016), 145701.

\bibitem{Ner81}
Nerman, O. On the convergence of supercritical general ({C}-{M}-{J}) branching
  processes. \emph{Z. Wahrsch. verw. Gebiete}, {\bf 57} (1981), 365--395.

\bibitem{Pasko}
Pasko, A. Yu., Likhachev, A.A., Koval Yu. N. and Kolomytsev, V.I.
2D Fourier Analysis and its Application to Study of Scaling Properties and Fractal Dimensions of 
$\varepsilon$-Martensite Distribution in γ-Matrix of Fe-Mn-Si Alloy,
\emph{ J. Phys. IV France}, {\bf 07} (1997),  C5-435--C5-440.
 
\bibitem{Planes13} 
Planes, A., Manosa L. and Vives, E. Acoustic emission in martensitic transformations \emph{J. Alloys and Compounds} {\bf 577S} (2013), S699--S704.

\bibitem{Salje} 
Salje, E.K.H., Koppensteiner, J., Reinecker, M., Schranz, W. and Planes, A. 
Jerky elasticity: Avalanches and the martensitic transition in $\textrm{Cu}_{74.08}\textrm{Al}_{23.13}\textrm{Be}_{2.79}$ shape-memory alloy,
\emph{Appl. Phys. Lett.} {\bf 95} (2009), 231908.

\bibitem{Song}
Song, Y., Chen, X., Dabade, V., Shield, T.W. and James, R.D. Enhanced reversibility and unusual microstructure 
in a phase-transforming material. \emph{Nature}, {\bf 502} (2013), 85--88.

\bibitem{Torrens16} 
Torrens, J., Illa, X., Vives, E. and Planes, A.
Geometrical model for martensitic phase transitions: understanding criticality and weak universality during 
microstructure growth, \emph{Phys. Rev. E}, {\bf 95} (2017), 013001.

\bibitem{Uch}
Uchiyama, K. Spatial growth of a branching process of particles living in $\mathbb{R}^d$. \emph{Ann. Probab.}
{\bf 10} (1982), 896--918.

\end{thebibliography}
\end{document}